\documentclass[11pt]{article}

    \usepackage{amsfonts}
    \usepackage{amsmath}
    \numberwithin{equation}{section}
    \usepackage{hyperref}
    \usepackage[capitalise]{cleveref}
    \usepackage{CJK}
    \usepackage{amsfonts,amssymb,amsmath,mathrsfs,multirow,booktabs}
    \usepackage{xcolor,soul}
\usepackage{color,latexsym,amsfonts}
 \newcommand{\qed}{\hfill\rule{2mm}{3mm}\vspace{4mm}}

 \newtheorem{theorem}{Theorem}[section]
 \newtheorem{lemma}[theorem]{Lemma}
 \newtheorem{corollary}[theorem]{Corollary}
 \newtheorem{proposition}[theorem]{Proposition}
 \newtheorem{example}[theorem]{Example}
 \newtheorem{Definition}[theorem]{Definition}
 \newtheorem{remark}[theorem]{Remark}
 \newtheorem{condition}[theorem]{Condition}
 \newtheorem{conjecture}[theorem]{Conjecture}

 \def\blemma{\begin{lemma}\sl{}\def\elemma{\end{lemma}}}
 \def\btheorem{\begin{theorem}\sl{}\def\etheorem{\end{theorem}}}

 \def\bproposition{\begin{proposition}\sl{}\def\eproposition{\end{proposition}}}
 
 \def\bcondition{\begin{condition}\sl{}\def\econdition{\end{condition}}}
 \def\bexample{\begin{example}\rm{}\def\eexample{\end{example}}}

 \def\beqlb{\begin{eqnarray}}\def\eeqlb{\end{eqnarray}}
 \def\beqnn{\begin{eqnarray*}}\def\eeqnn{\end{eqnarray*}}

 \def\mbb{\mathbb}\def\mbf{\mathbf}

 \def\<{\langle}\def\>{\rangle}

 \def\ar{&\!\!}

 \def\eqref#1{{\rm(\ref{#1})}}

 \def\proof{\noindent{\it
 Proof.~}}\def\qed{\hfill$\Box$\medskip}

\def\e{{\mbox{\rm e}}}
\def\<{\left<}\def\>{\right>}

 \def\mbb{\mathbb}
  \def\mbf{\mathbf}
\newcommand{\dd}{\mathrm{d}}

\newfam\msbmfam\font\tenmsbm=msbm10\textfont
\msbmfam=\tenmsbm\font\sevenmsbm=msbm7
\scriptfont\msbmfam=\sevenmsbm

\def\<{\left<}\def\>{\right>}
\def\({\left(}\def\){\right)}

\begin{document}

\centerline{\Large\bf Extinction, explosion and contraction}
	
	\smallskip
	
	\centerline{\Large\bf for time-inhomogeneous SDEs with jumps}
	
	\bigskip

\centerline{Shukai Chen$^1$, Xu Yang$^2$ and  Xiaowen Zhou$^3$}

\medskip
	
	{\small\it

       \centerline{$^1$School of Mathematics and Statistics, Fujian Normal University, Fuzhou, China.}

       \centerline {\tt skchen@fjnu.edu.cn}

		\smallskip
        \centerline{$^2$School of Mathematics and Information Science,
			North Minzu University, Yinchuan, China. }
		
		\centerline{{\tt xuyang@mail.bnu.edu.cn}}
		
		\smallskip
        \centerline{$^3$Department of Mathematics and Statistics, Concordia
			University, Montreal, Canada.}

        \centerline{\tt xiaowen.zhou@concordia.ca}
		
	}

	\bigskip
	
	{\narrower{\narrower
			
			\noindent{\textit{Abstract:}}
			 For a class of time-inhomogeneous SDEs with jumps, we establish criteria for the existence and uniqueness of the nonnegative solutions, and examine the  extinction, the explosion  together with the contractivity of the solutions, which generalize and improve upon earlier results  in the literature. As an application, we study the aforementioned properties for a class of mean field SDEs.
			
			
			\smallskip
			
			\noindent{\textit{MSC {\rm(2020)} Subject Classification:}} 60H20, 60J50, 60J25.
			
			\par}\par}

	\bigskip

\section{Introduction}\label{intro11}
\setcounter{equation}{0}

\subsection{Background}
The theory of stochastic integrals and stochastic differential equations (SDEs for short) was pioneered by It\^{o} in 1942. It was first applied to solve a fundamental problem posed by Kolmogorov on determining Markov processes in 1931. Compared with ordinary differential equations, SDEs are designed to model differential equations perturbed by random noise. Today, this theory plays an increasingly vital role in diverse fields such as dynamical systems, mathematical finance, and biology. Initially, researchers primarily focused on stochastic equations and stochastic integrals driven by Brownian motion. Let  $(B_t)_{t\geq 0}$ be a standard Brownian motion. Given  functions $b$ and $\sigma$ on $[0,\infty)\times[-\infty,\infty]^d$ for $d\ge1$, the  SDE is of the form
\beqlb\label{SDE brownian}
\dd X_t=b(t,X_t)\,\dd t+\sigma(t,X_t)\,\dd B_t,
\eeqlb
or equivalently,
\beqnn
X_t=X_0+\int_0^tb(s,X_s)\,\dd s+\int_0^t\sigma(s,X_s)\,\dd B_s.
\eeqnn
Then a natural question  is the existence and uniqueness of its solutions. For foundational contributions to \eqref{SDE brownian}, the literature is extensive; we refer to  \cite{KS87, LeG, Mao08, O03, RY91, RW94, SV79, YW71} for further details. Since Brownian motion is a special case of L\'evy processes, subsequent research has extended to stochastic equations driven by general L\'evy processes $(L_t)_{t\ge0}$, giving rise to SDEs of the form
\beqlb\label{SDE jumps}
\dd X_t=b(t,X_t)\,\dd t+\sigma(t,X_t)\,\dd L_t,
\eeqlb
which are often referred to as stochastic equations with jumps. By the L\'evy-It\^{o} decomposition, \eqref{SDE jumps} can be seen as SDEs driven by certain Poisson random measures. References for this line of research include, among others, \cite{B03, IW89, JS03, P05, Si05}. We say that the SDEs \eqref{SDE brownian} and \eqref{SDE jumps} are time-homogeneous if the coefficients $b$ and $\sigma$ do not depend explicitly on the time parameter $t$. Otherwise, it is termed a time-inhomogeneous SDE.

In addition to the existence and uniqueness of solutions, other theoretical properties of the solution as a stochastic process have also attracted considerable attention, particularly among researchers in mathematical  biology. These properties include its long-term behavior, such as whether the solution diverges to infinity, hits zero, or exhibits ergodicity. For this purpose, it is common practice to consider stochastic equations with nonnegative solutions. Among others, we refer to \cite{BHS08, FZ23, FJKR23, HN18, HNC21, HNTU25,  LZ21, LYZh, MZ24, M11, NY17,RXYZ19, ZY09}  for associated results of some stochastic models for population dynamics.

In this work, we focus on a class of time-inhomogeneous processes
 as  nonnegative solutions to time-inhomogeneous stochastic equations
and investigate three long time behaviors of such processes for which
the known  results  are relatively limited.
Such processes are formally defined as follows. Let $(\Omega, \mathcal{F}, (\mathcal{F}_t)_{t\geq 0}, \mathbf{P})$ be a filtered probability space satisfying the usual conditions, and let $\mathbf{E}$ denote the expectation with respect to $\mathbf{P}$.  Let $\{W(\dd s,\dd u): s,u\ge0\}$ be a time-space $(\mathcal{F}_t)$-Gaussian white noise with density $2b_1(s-)\dd s\dd u$ and  let $\{N(\dd s,\dd z,\dd u): s,z,u>0\}$ be a time-space $(\mathcal{F}_t)$-Poisson
random measure with intensity $b_2(s-)\dd s\mu(\dd z)\dd u$, where $\mu$ is a $\sigma$-finite measure on $(0,\infty)$ satisfying $\int_0^\infty (z\wedge z^2)\mu(\dd z)<\infty$, and $b_1, b_2$ are some positive c\`{a}dl\`{a}g Borel functions. Denote by  $\tilde{N}(\dd s,\dd z,\dd u)$ the compensated measure of $N(\dd s,\dd z,\dd u)$.  We further assume that $W$ and $N$ are independent. We consider
\begin{equation}\label{ii1.5}
X_t =X_0+\int_0^t\gamma_0({s-},X_{s-})\dd s
+\int_0^t\int_0^{\gamma_1(X_{s-})}W(\dd s,\dd u)
+\int_0^t\int_0^\infty\int_0^{\gamma_2(X_{s-})}z\tilde{N}(\dd s,\dd z,\dd u),
\end{equation}
where $X_0>0$,
$\gamma_0$ and $\gamma_1,\gamma_2\ge0$ are Borel functions on $[0,\infty)\times[0,\infty)$ and $[0,\infty)$, respectively. Here and in what follows, we make the convention that
$$
\int_x^y=\int_{(y,x]} ~~\text{and}~~ \int_x^\infty=\int_{(x,\infty)},\quad x\ge y\in\mbb{R}.
$$
We refer to \cite{Wal86} for the stochastic integral with respect to white noise and to \cite{IW89} for the stochastic integral with respect to Poisson random measures. One could replace the diffusion term in \eqref{ii1.5} by the stochastic integral $\int_0^t\sqrt{2\gamma_1(X_{s-})b_1(s-)}\,\dd B_s$ using a one-dimensional $(\mathcal{F}_t)$-Brownian motion $(B_t)_{t\ge0}$. The resulting equation defines an equivalent process for any fixed $X_0>0$. A key advantage of using white noise in \eqref{ii1.5} over Brownian motion is the former's temporal independence. As shown in \cite{DL12},  this property naturally ensures the equivalence of the flows generated by the solution.

We emphasize that \eqref{ii1.5} features several interesting  stochastic population models such as branching processes. \cite{DL06, DL12, FL10, LM11,LP12} studied the pathwise uniqueness and strong solutions for a class of time-homogeneous SDEs
driven by spectrally positive L\'evy noises with nonnegative solutions. The origin of these processes (solutions to stochastic equations) can be traced back to branching processes, which were initially developed to model population extinction. In particular, the construction of continuous state and continuous time branching processes (CB-processes for short) through stochastic equations can be viewed as the unique strong solution to the Dawson-Li stochastic equation:
\begin{equation}\label{CB-process}
X_t =X_0-\int_0^tbX_{s-}\dd s
+\int_0^t\int_0^{X_{s-}}B(\dd s,\dd u)
+\int_0^t\int_0^\infty\int_0^{X_{s-}}z\tilde{M}(\dd s,\dd z,\dd u),
\end{equation}
where $b\in\mathbb{R}, B$ denotes a Gaussian white noise and $\tilde{M}$ denotes an independent
compensated Poisson
random measure. Later, \cite{FL22} derived the stochastic equation
\begin{equation}\label{CBVE-process}
X_t =X_0-\int_0^tX_{s-}c_0(\dd s)
+\int_0^t\int_0^{X_{s-}}B_1(\dd s,\dd u)
+\int_0^t\int_0^\infty\int_0^{X_{s-}}z\tilde{M}_1(\dd s,\dd z,\dd u),
\end{equation}
for the CB-process in varying environments, where $B_1$ is a Gaussian white noise with density $2c_1(\dd s)\dd u$, and $\tilde{M}_1$ is an independent
compensated Poisson
random measure with intensity $m(\dd s,\dd z)\dd u$. A deep understanding of such
models in varying environments is important, as they provide the basis for further study of models in random
environments. We refer to \cite{BCM19, HLX18, PP17} and references therein on this research program. Meanwhile, research on time-inhomogeneous SDEs has found broad applications in fields such as machine learning and image recognition. \cite{LGZSS23a, LGZSS23b} leveraged time-inhomogeneous SDEs to develop the so-called IR-SDE method, which enables the transformation of low-quality images into high-quality ones. \cite{ZSHLCLW24} introduced the SDE-EDG method based on time-inhomogeneous SDEs to tackle the challenge of time-varying distribution shifts in machine learning applications. Additionally, \cite{KSZ20} proposed a novel uncertainty-aware neural network from the perspective of  stochastic dynamical systems. Their goal was to provide efficient uncertainty quantification methods for deep learning. By interpreting the forward propagation of neural networks as state evolution within a dynamical systems framework, their approach incorporates Brownian motion to capture model uncertainties. A class of time-inhomogeneous SDEs  also plays a central role in this model.

\subsection{Outline of main results}

This paper is structured around the following three main results.

Our first main result, Theorem \ref{strong solution of X}, establishes the existence and pathwise uniqueness of strong solutions to \eqref{ii1.5}. Consequently, it  follows from \cite[Theorems 4.1, 4.2, and 7.1]{EK86} that the solution is a Markov process.
The proof of pathwise uniqueness is a modification
of the classical Yamada-Watanabe argument.
The existence is established by considering piecewise time-homogeneous SDEs.

The second main result, Theorems \ref{t1.3}--\ref{t1.4c},
gives rather sharp conditions on extinction/non-extinction and explosion/non-explosion,
and generalize the corresponding conditions  in \cite{LYZh} for  time-homogeneous SDEs  to  time-inhomogeneous SDEs.
The proofs are inspired by Chen's criteria for Markov jump processes
(see, e.g., \cite{chen04,MT93}).
The key to our approach is to find an
appropriate test function that yields  the best possible result. However, these functions may not always have an obvious intuitive interpretation and are often chosen in an ad hoc manner.
The test functions we selected in this paper are often logarithm type functions, power functions, exponential functions
 or their composition.
In particular, the test functions for extinction
and explosion are composition of logarithm type functions
and exponential functions, which are different  from those in previous papers (see Lemma \ref{t1.7} and the proof of
Theorem \ref{t1.4}).

 Our third main result,
Theorem \ref{exponential ergod}, establishes the contractivity for the unique solution $X$ to \eqref{ii1.5} and generalizes the results of \cite{LLWZ25}. Our conditions are more general than those in \cite{LLWZ25}, which necessitates a careful treatment of the effects of time-inhomogeneity and nonlinearity.

Finally, to better illustrate the above results, we provide an SDE in Example 4.4 that has a unique solution and the solution is non-explosive and contractive.
As an important application, we also provide an example of a class of killed mean field dependent SDEs in Example \ref{e4.2}.
A systematic discussion of killed distribution-dependent SDEs can be found in \cite[Chapter 7]{WR25}.  Example \ref{e4.2} shows that  the results of this paper can be used to show that such models can be regarded as time-inhomogeneous Markov processes, which are non-explosive, non-extinctive, and possess exponential contractivity properties.

The remaining of the paper is organized as follows. In Section 2, we exploit the existence and uniqueness of solutions to \eqref{ii1.5}. The extinction and explosion behaviors are established in Section 3. The contractive properties of the associated transition semigroups are proved in Section 4.

{\bf Notation.}
Let $X$ be the solution to equation \eqref{ii1.5}. For $0 < a < b $, define the stopping times
\beqlb\label{3.2}
\tau_b^+ := \inf\{ t \ge 0 : X_t \ge b \}, \qquad \tau_a^- := \inf\{ t \ge 0 : X_t \le a \}.
\eeqlb
We then define the limiting stopping times as
$
\tau_0 := \lim_{a \to 0} \tau_a^-$ and $\tau_\infty := \lim_{b \to \infty} \tau_b^+.
$ For $m, n \ge 1 $, we define the stopping time
$
\tau_{m,n} := \tau^-_{1/m} \wedge \tau^+_n.
$ Separately, for a nonnegative Borel function $g$, we define
$$
g((0, t]) := \int_0^t g(s) \, \dd s, \qquad g((0, \infty)) := \int_0^\infty g(s) \, \dd s.
$$

Let $ \mbb{R}_+ = [0, \infty) $ and $ \mbb{R} = (-\infty, \infty)$. For a set $ A \subset \mbb{R}^d $ ($d \ge 1$), we denote by  $ C^2(A) $ the space of real-valued, twice continuously differentiable functions on $A$. For the process $X$ starting from $ x \in \mbb{R}^d $, we denote its law by $ \mathbf{P}^x(\cdot)$ and the corresponding expectation by $\mathbf{E}^x(\cdot) $. For a differentiable function $f$ on $\mbb{R}^d$, define
$$
\Delta_z f(x) := f(x+z) - f(x), \quad D_z f(x) := \Delta_z f(x) - z \cdot \nabla f(x),
$$
where $ \nabla f$ denotes the gradient of $f$, and  $\cdot$ is the dot product in $\mbb{R}^d$. In the one-dimensional case ($d=1$), this definitions reduce to
$$
\Delta_z f(x) = f(x+z) - f(x), \quad D_z f(x) = \Delta_z f(x) - z f'(x).
$$

The generator $ (\mathcal{L}_s)_{s \ge 0} $ of the one-dimensional process $X$ is defined as follows. For any $s \ge 0$ and $f \in C^2(\mbb{R}_+)$,
\beqlb\label{generator of X}
\mathcal{L}_s f(x)
\ar=\ar \gamma_0(s, x) f'(x) + \gamma_1(x) b_1(s) f''(x) + \gamma_2(x) b_2(s) \int_0^\infty D_z f(x) \mu(\dd z) \cr
\ar=\ar \gamma_0(s, x) f'(x) + \gamma_1(x) b_1(s) f''(x)\cr
\ar\ar+ \gamma_2(x) b_2(s) \int_0^\infty z^2 \mu(\dd z) \int_0^1 f''(x + z u)(1 - u) \dd u,
\eeqlb
where the second equality follows from a Taylor expansion.

\section{Existence and uniqueness}

We first introduce the following two definitions.

{\begin{Definition}\label{Def of weak solu}
By a solution to SDE \eqref{ii1.5}, we mean a c\`{a}dl\`{a}g process $X$ satisfying \eqref{ii1.5} up to time $\tau_n := \tau^-_{1/n}\wedge \tau^+_n$ for each $n \ge 1$ and $X_t = \limsup_{n\rightarrow\infty} X_{\tau_n^-}$ for
$t \ge \tau := \lim_{n\rightarrow\infty} \tau_n$.
\end{Definition}
Then, by Definition \ref{Def of weak solu}, both boundary points $0$ and $\infty$ are absorbing for $X$ .

\begin{Definition}(Pathwise uniqueness)
The solution of SDE \eqref{ii1.5} is said to be pathwise unique if, for any two solutions $X$ and $\tilde{X}$ satisfying \eqref{ii1.5} with the same initial value $X_0=\tilde{X}_0$, it holds that
$$
\mathbf{P}\left(X_t=\tilde{X}_t,~~\text{for~all}~ t\in[0,\tau\wedge\tilde{\tau}]\right)=1,
$$
where $\tau$ is the stopping time defined as in Definition \ref{Def of weak solu} for $X$, and $\tilde{\tau}$ is defined analogously for $\tilde{X}$.
\end{Definition}

Now we give sufficient conditions on the functions $\gamma_i, i=0,1,2$ under which SDE \eqref{ii1.5} has a pathwise unique solution $X$.

\bcondition\label{cond for solution}

(i)~The function $x\mapsto\gamma_0(\cdot,x)$ is locally Lipschitz in the sense that there is a positive and measurable function $b_0(s)$ and for each closed interval $[a,b]$ with $0<a<b<\infty$, there is a constant $c_{a,b}>0$ such that for any $x,y\in [a,b]$ and $s>0$,
$$
|\gamma_0(s,x)-\gamma_0(s,y)|\le c_{a,b}|x-y|b_0(s),
$$
and for any $t>0$ and $x_0>0$,
\beqnn
\sup_{0<s\le t,0\le x\le x_0} |\gamma_0(s,x)|<\infty,\quad
\sup_{0<s\le t}[b_0(s)+b_1(s)+b_2(s)]<\infty.
\eeqnn

(ii)~The functions $x\mapsto\gamma_1(x)$ and $x\mapsto\gamma_2(x)$ are locally Lipschitz; that is, for each closed interval $[a,b]$ with $0<a<b<\infty$, there is a constant $c_{a,b}>0$ such that for any $x,y\in [a,b]$,
$$
|\gamma_1(x)-\gamma_1(y)|+|\gamma_2(x)-\gamma_2(y)|\le c_{a,b}|x-y|.
$$
Moreover, $x\mapsto\gamma_2(x)$ is increasing, and for any $x_1>x_2>0$, $\sup_{x_1\le x\le x_2} [\gamma_1(x)+\gamma_2(x)]<\infty$.

\econdition

\btheorem\label{strong solution of X}
Suppose that Condition \ref{cond for solution} holds. Then for any $X_0 = x > 0,$ there exists a pathwise unique solution $X$ to SDE \eqref{ii1.5}.
\etheorem

With Theorem \ref{strong solution of X} established, the Markov property of the solution $X$ follows by applying \cite[Theorems 4.1, 4.2 and 7.1]{EK86}.

Before presenting the detailed proof, we briefly outline the main strategy. The proof proceeds in two steps: establishing the existence of a solution (Proposition \ref{ii1.8}) and proving the pathwise uniqueness of the solution (Proposition \ref{ii1.9}).

{\bf Step 1: Existence of a solution.} First, we construct a family of discrete approximation processes by discretizing the functions $\gamma_0,\gamma_1,\gamma_2, b_1$ and $b_2$. Then on each subinterval, the corresponding SDE becomes a time-homogeneous SDE. This allows us to apply  classical results to ensure the existence of a nonnegative strong solution on every small interval. These local solutions are then pieced together to obtain a global discretized solution. Next, we verify the tightness of this family of processes. The   proof  here primarily follows the approach in \cite{XY19}. With tightness established, we employ tools such as Skorokhod's representation theorem to prove the convergence of a subsequence if necessary, ensuring that the limiting process satisfies \eqref{ii1.5}.

{\bf Step 2:  Pathwise uniqueness.} The proof of pathwise uniqueness follows a standard argument, similar to classical methods in the literature, see \cite{DL12,FL10} for instance.

We first introduce the following notation. For any $m,n\ge1$, we define the following sequences,
\beqlb\label{discrete sequence}
\ar\ar
b_i^n(s)=\sum_{k=1}^\infty\mathbf{1}_{(s^n_{k-1},s^n_k]}(s)b_i(s^n_{k-1}),\quad i=1,2,\cr
\ar\ar
\gamma^m_i(x)=\gamma_i(m)\mathbf{1}_{\{m<x<\infty\}}+\gamma_i(x)\mathbf{1}_{\{1/m\le x\le m\}}+\gamma_i(1/m)\mathbf{1}_{\{0\le x<1/m\}},\quad i=1,2,\cr
\ar\ar \gamma_0^{n,m}(s,x)=\sum_{k=1}^\infty\mathbf{1}_{(s^n_{k-1},s^n_k]}(s)
\gamma^m_0\left(s^n_{k-1},x\right),
\eeqlb
where $s^n_k=k/n$ and for all $s\ge0$,
$$
\gamma^m_0(s,x)=\gamma_0(s,m)\mathbf{1}_{\{m<x<\infty\}}+\gamma_0(s,x)\mathbf{1}_{\{1/m\le x\le m\}}+\gamma_0(s,1/m)\mathbf{1}_{\{0\le x<1/m\}}.
$$
Then under Condition \ref{cond for solution} (ii), for each $m\ge1$, there is a constant $C>0$ such that
 \beqlb\label{5.1}
\gamma^m_1(x)+\gamma^m_2(x)\le C,\qquad x>0
 \eeqlb
and
 \beqlb\label{5.2}
|\gamma^{n,m}_0(s,x)|
 \ar\le\ar
\sum_{k=1}^\infty\mathbf{1}_{(s^n_{k-1},s^n_k]}(s)
\Big|\gamma^m_0\left(s^n_{k-1},x\right)-\gamma_0\left(s^n_{k-1},1\right)
+\gamma_0\left(s^n_{k-1},1\right)\Big| \cr
\ar\le\ar
\sum_{k=1}^\infty\mathbf{1}_{(s^n_{k-1},s^n_k]}(s)
\big[C|b_0(s^n_{k-1})|
+|\gamma_0(s^n_{k-1},1)|\big],\quad x>0.
 \eeqlb

Let $W_n(\dd s,\dd u)$
be a Gaussian white noise with density $2b^n_1(s)\dd s\dd u$, and let
$N_n(\dd s,\dd z,\dd u)$ be a Poisson
random measure with intensity $b^n_2(s)\dd s\mu(\dd z)\dd u$.
For any $n,m\ge1$, consider the following SDE:
\beqlb\label{ii1.6}
X_t \ar=\ar
X_0+\int_0^t\gamma^{n,m}_0(s,X_{s-})\dd s
+\int_0^t\int_0^{\gamma^m_1(X_{s-})}W_n(\dd s,\dd u)\cr
\ar\ar
+\int_0^t\int_0^\infty\int_0^{\gamma^m_2(X_{s-})}z\tilde{N}_n(\dd s,\dd z,\dd u),
\eeqlb
where $\tilde{N}_n(\dd s,\dd z,\dd u)=N_n(\dd s,\dd z,\dd u)-b^n_2(s)\dd s\mu(\dd z)\dd u$.

\blemma\label{ii1.7}
Suppose that Condition \ref{cond for solution} holds. Then the following assertions hold:

(a) For any $n,m\ge1$, there is a unique nonnegative strong solution $(X^{n,m}_t)_{t\ge0}$ to \eqref{ii1.6}.

(b) Fix $m\ge1$. The family $\{(X^{n,m}_t)_{t\ge0}: n\ge1\}$ is tight in the Skorokhod space $D([0,\infty),\mathbb{R}_+)$.
\elemma

\proof
We first prove (a). Note that on each interval $(s^n_{k-1},s^n_k]$, \eqref{ii1.6} reduces to a time-homogeneous SDE.  By \cite[Theorem 2.5]{DL12}, there exists a unique nonnegative strong solution $(X^{n,m}_t)_{t\in(s^n_{k-1},s^n_k]}$. By piecing these solutions together, we obtain a unique nonnegative strong solution $(X^{n,m}_t)_{t\ge0}$ to the SDE \eqref{ii1.6}.

We now prove part (b) using \cite[Lemma 2.2]{XY19}. Specifically, we verify the following two conditions:

\noindent(i) For every $\eta>0, m\ge1$ and $T>0$, there is a constant $\lambda_{\eta,T,m}>0$ such that
\beqlb\label{YX cond1}
\inf_{n\ge1}\mathbf{P}
\left\{\sup_{t\in[0,T]}X^{n,m}_t\le \lambda_{\eta,T,m}\right\}\ge 1-\eta;
\eeqlb
(ii) Let $f$ be in a subalgebra that is dense in $C_b(\mathbb{R}_+)$ (the bounded continuous functions space on $\mathbb{R}_+$), there is a process $(g^{n,m}_t)_{t\ge0}$ so that
\beqlb\label{YX cond2.1}
f(X^{n,m}_t)-\int_0^tg^{n,m}_s\,\dd s
\eeqlb
is a martingale and
\beqlb\label{YX cond2.2}
\sup_{t\in[0,T]}\mathbf{E}\left\{|f(X^{n,m}_t)|+|g^{n,m}_t|\right\}<\infty,\quad
\sup_{n\ge1}\mathbf{E}\left\{\left[\int_0^T|g^{n,m}_t|^p\,\dd t\right]^{1/p}\right\}<\infty
\eeqlb
for each $n\ge1, T>0$ and some $p>1$.

By \eqref{ii1.6}, we have
\beqlb\label{est of tig}
\mathbf{E}^x\left[\sup_{0\le t\le T}X^{n,m}_t\right]
\ar\le\ar x+\mathbf{E}^x\left[\int_0^T|\gamma_0^{n,m}(s,X^{n,m}_{s-})|\,\dd s\right]
+\mathbf{E}^x\left[\sup_{0\le t\le T}\left|\int_0^t\int_0^{\gamma_1^{m}(X^{n,m}_{s-})}W_n(\dd s,\dd u)\right|\right]\cr
\ar\ar
+\mathbf{E}^x\left[\sup_{0\le t\le T}\left|\int_0^t\int_0^1\int_0^{\gamma_2^m(X^{n,m}_{s-})}z\tilde{N}_n(\dd s,\dd z,\dd u)\right|\right]\cr
\ar\ar
+\mathbf{E}^x\left[\sup_{0\le t\le T}\left|\int_0^t\int_1^\infty\int_0^{\gamma_2^m(X^{n,m}_{s-})}z\tilde{N}_n(\dd s,\dd z,\dd u)\right|\right].
\eeqlb
By Doob's inequality and \eqref{5.1},
\beqlb\label{est of expect gauss}
\mathbf{E}^x\left[\sup_{0\le t\le T}\left|\int_0^t\int_0^{\gamma_1^{m}(X^{n,m}_{s-})}W_n(\dd s,\dd u)\right|^2\right]
\ar\le\ar 4\mathbf{E}^x\left[\int_0^Tb_1^n(s)\gamma_1^m(X^{n,m}_s)\dd s\right]\cr
\ar\le\ar
4C\int_0^Tb_1^n(s)\dd s
\eeqlb
and
\beqlb\label{est of expect jump1}
\mathbf{E}^x\left[\sup_{0\le t\le T}\left|\int_0^t\int_0^1\int_0^{\gamma_2^m(X^{n,m}_{s-})}z\tilde{N}_n(\dd s,\dd z,\dd u)\right|^2\right]
\ar\le\ar 4\int_0^1z^2\mu(\dd z)\mathbf{E}\left[\int_0^Tb_2^n(s)\gamma_2^m(X^{n,m}_{s-})\dd s\right]\cr
\ar\le\ar
4C\int_0^1z^2\mu(\dd z)\int_0^Tb_2^n(s)\dd s.
\eeqlb
For $z>1$, by using the fact that $\tilde{N}_n(\dd s,\dd z,\dd u)=N_n(\dd s,\dd z,\dd u)-b^n_2(s)\dd s\mu(\dd z)\dd u$, we have
\beqlb\label{est of expect jump2}
\ar\ar\mathbf{E}^x\left[\sup_{0\le t\le T}\left|\int_0^t\int_1^\infty\int_0^{\gamma_2^m(X^{n,m}_{s-})}z\tilde{N}_n(\dd s,\dd z,\dd u)\right|\right]\cr
\ar\ar\le \mathbf{E}^x\left[\int_0^T\int_1^\infty\int_0^{\gamma_2^m(X^{n,m}_{s-})}z\,N_n(\dd s,\dd z,\dd u)\right]
+\mathbf{E}^x\left[\int_0^T\int_1^\infty\int_0^{\gamma_2^m(X^{n,m}_{s-})}z\,b^n_2(s)\dd s\mu(\dd z)\dd u\right]\cr
\ar\ar
\le 2C\int_1^\infty z\,\mu(\dd z)\int_0^Tb_2^n(s)\,\dd s.
\eeqlb
Collecting estimates \eqref{5.2}, \eqref{est of expect gauss}, \eqref{est of expect jump1} and \eqref{est of expect jump2} with Condition \ref{cond for solution} (i) and \eqref{est of tig}, we obtain
$$
\lambda_{m,T}=\sup_{n\ge1}\mathbf{E}^x\left[\sup_{0\le t\le T}X^{n,m}_t\right]<\infty
$$
for any $m\ge1$ and $T>0$. Set $
\lambda_{\eta,T,m}:=\lambda_{m,T}\eta^{-1}.$     Then \eqref{YX cond1} follows from Markov inequality.
We now verify \eqref{YX cond2.1} and \eqref{YX cond2.2}. By It\^{o}'s formula and \eqref{ii1.6}, for each $F\in C^2(\mathbb{R}_+)$ with compact supports,
\beqnn
F(X^{n,m}_t)\ar=
\ar F(X_0)+\int_0^tF'(X^{n,m}_{s-})\gamma_0^{n,m}(s,X^{n,m}_{s-})\,\dd s+\int_0^tF''(X^{n,m}_{s-})\gamma^m_1(X^{n,m}_{s-})b^{n}_1(s)\,\dd s \cr
\ar\ar +\int_0^t\int_0^\infty D_zF(X^{n,m}_{s-})\cdot b^{n}_2(s)\gamma_2^m(X^{n,m}_{s-})\,\dd s\mu(\dd z)+\mbox{martingale},
\eeqnn
with the operator $D_z$ defined at the end of Section 1. By Taylor expansion and $F\in C^2(\mathbb{R}_+)$ with compact supports, for any $x,z\ge0$, we have
$$
\int_0^1|D_zF(x)|\,\mu(\dd z)=\frac{1}{2}\int_0^1z^2F''(\zeta)\,\mu(\dd z)<\infty
$$
for $\zeta\in(x,x+z)$ and
$$
\int_1^\infty|D_zF(x)|\,\mu(\dd z)\le \left(2\sup_x|F(x)|+\sup_x|F'(x)|\right)\int_1^\infty z\mu(\dd z)<\infty.
$$
This completes the proof.\qed

\bproposition\label{ii1.8}
Suppose that Condition \ref{cond for solution} holds. Then there is a solution $(X_t)_{t\ge0}$ to the SDE \eqref{ii1.5}.
\eproposition
\proof  By Lemma \ref{ii1.7} (b), $(X^{n,m}, W_n, N_n)$ in \eqref{ii1.6} converges in distribution to some limit $(X^m, W, N)$  as $n\rightarrow\infty$ (along a subsequence if necessary). By Skorokhod's representation theorem, there exist a filtered probability space satisfying the usual conditions, random vectors $(\hat{X}^m,\hat{W},\hat{N})$ and $\{(\hat{X}^{n,m}, \hat{W}_n, \hat{N}_n): n \ge 1\}$ defined on this space with the same distribution as $(X^m, W, N)$ and $\{(X^{n,m}, W_n, N_n): n \ge 1\}$, respectively. Moreover, $(\hat{X}^{n,m}, \hat{W}_n, \hat{N}_n)$ converges to $(\hat{X}^m,\hat{W},\hat{N})$ almost surely as $n\rightarrow\infty$.

By Lemma \ref{ii1.7}   (a), \eqref{ii1.6} admits a strong solution. Consequently, $(\hat{X}^{n,m}, \hat{W}_n, \hat{N}_n)$ satisfies the following equation:
\beqlb\label{ii1.10}
\hat{X}^{n,m}_t \ar=\ar X_0 + \int_0^t\gamma^{n,m}_0(s-,\hat{X}^{n,m}_{s-})\dd s
+\int_0^t \int_0^{\gamma^m_1(\hat{X}^{n,m}_{s-})} \hat{W}_n(\dd s,\dd u)\cr
\ar\ar\,
 + \int_0^t\int_0^\infty\int_0^{\gamma^m_2(\hat{X}^{n,m}_{s-})}z\tilde{\hat{N}}_n(\dd s,\dd z,\dd u),
\eeqlb
where $\tilde{\hat{N}}_n(\dd s,\dd z,\dd u):=\hat{N}_n(\dd s,\dd z,\dd u)-b^n_2(s)\dd s\mu(\dd z)\dd u$.
As established in \eqref{est of expect gauss}, for any $T>0,m\ge1$ and $\delta>0$, Markov inequality yields
\beqnn
\sup_{n\ge1}\mathbf{P}\left(\sup_{0\le t\le T}\left[\left|\int_0^t
 \int_0^{\gamma^m_1(\hat{X}^{n,m}_{s-})} \hat{W}_n(\dd s,\dd u)\right|\right]\ge C_{T,\delta,m}\right)\le\delta,
\eeqnn
where
$$
C_{T,\delta,m}=\sup_{n\ge1}\mathbf{E}\left(\sup_{0\le t\le T}\left|\int_0^t
 \int_0^{\gamma^m_1(\hat{X}^{n,m}_{s-})} \hat{W}_n(\dd s,\dd u)\right|\right)\delta^{-1}.
$$
Similarly, by \eqref{est of expect jump1} and \eqref{est of expect jump2}, we have
\beqnn
\sup_{n\ge1}\mathbf{P}\left(\sup_{0\le t\le T}\left[\left|\int_0^t\int_0^\infty\int_0^{\gamma^m_2(\hat{X}^{n,m}_{s-})}z\tilde{\hat{N}}_n(\dd s,\dd z,\dd u)\right|\right]\ge C_{T,\delta,m}\right)\le\delta
\eeqnn
with
$$
C_{T,\delta,m}=\sup_{n\ge1}\mathbf{E}\left(\sup_{0\le t\le T}\left|\int_0^t\int_0^\infty\int_0^{\gamma^m_2(\hat{X}^{n,m}_{s-})}z\tilde{\hat{N}}_n(\dd s,\dd z,\dd u)\right|\right)\delta^{-1}.
$$
By \cite[Lemma 2.4]{XY19}, for any $t>0$, the following convergences hold in probability as $n\rightarrow\infty$:
$$
\int_0^t \int_0^{\gamma^m_1(\hat{X}^{n,m}_{s-})} \hat{W}_n(\dd s,\dd u)
\longrightarrow\int_0^t\int_0^{\gamma^m_1(\hat{X}^m_{s-})}\hat{W}(\dd s,\dd u)
$$
and
$$
\int_0^t\int_0^\infty\int_0^{\gamma^m_2(\hat{X}^{n,m}_{s-})}z\tilde{\hat{N}}_n(\dd s,\dd z,\dd u)
\longrightarrow\int_0^t\int_0^\infty\int_0^{\gamma^m_2(\hat{X}^m_{s-})}z\tilde{\hat{N}}(\dd s,\dd z,\dd u).
$$
Moreover, by the dominated convergence theorem, the continuity of $x\mapsto\gamma_0^{n,m}(\cdot,x)$ and Condition \ref{cond for solution} (i),  we have
\beqnn
\int_0^t\gamma_0^{n,m}(s-,\hat{X}^{n,m}_{s-})\,\dd s\longrightarrow\int_0^t\gamma_0^{m}(s-,\hat{X}^m_{s-})\,\dd s
\eeqnn
almost surely as $n\rightarrow\infty$. Passing to the limit in \eqref{ii1.10} along an appropriate subsequence, we find that $(\hat{X}^m_t)_{t\ge0}$ satisfies
\beqnn
\hat{X}^m_t = X_0 + \int_0^t\gamma^{m}_0(s-,\hat{X}^m_{s-})\dd s
+\int_0^t\int_0^{\gamma^m_1(\hat{X}^m_{s-})}\hat{W}(\dd s,\dd u)+\int_0^t\int_0^\infty\int_0^{\gamma^m_2(\hat{X}^m_{s-})}z\tilde{\hat{N}}(\dd s,\dd z,\dd u).
\eeqnn
Finally, employing arguments analogous to those in the proof of \cite[Theorem 3.1]{LYZh}, we conclude that \eqref{ii1.5}  admits a solution $(X_t)_{t\ge0}$. \qed

\bproposition\label{ii1.9}
Suppose that Condition \ref{cond for solution} holds. Then the pathwise uniqueness of solutions holds for \eqref{ii1.5}.
\eproposition

\proof See the Appendix.\qed

\section{Extinction and explosion}\label{intro12}

In this section, we study conditions for
non-extinction/extinction and non-explosion/explosion
conditions of solutions to \eqref{ii1.5}, and we assume that the following assumptions hold:
\beqlb\label{3.0}
\sup_{x_1\le x\le x_2,\,0<s\le t}
[b_1(s)+b_2(s)+\gamma_1(x)+\gamma_2(x)+|\gamma_0(s,x)|]<\infty
 \eeqlb
for all $t>0$ and $x_2\ge x_1>0$.
The key to the proofs is the estimation of $\mathcal{L}_s$ given by \eqref{generator of X}.

\subsection{Main results}

Before presenting the main results, we introduce the following notation.
For $s,x>0$, define
 \beqlb\label{3.0a}
H(s,x)=
b_1(s)x^{-2}\gamma_1(x)
+b_2(s)\gamma_2(x)
\int_0^\infty z^2\mu(\dd z)\int_0^1(x+zu)^{-2}(1-u)\dd u.
 \eeqlb
We now state conditions for the non-extinction and non-explosion of solutions to \eqref{ii1.5}, which apply near zero and near infinity, respectively.

\btheorem\label{t1.3}
\begin{itemize}
\item[{\normalfont(i)}]
If for all $t>0$, there is a constant $c_0=c_0(t)\in(0,1)$ such that
 \beqlb\label{3.aaa}
\sup_{s\in(0,t]\,0<x<c_0}[\ln x^{-1}]^{-1}[-x^{-1}\gamma_0(s,x)+H(s,x)]
<\infty,
 \eeqlb
then $\mbf{P}\{\tau_0<\infty\}=0$.
\item[{\normalfont(ii)}]
If for all $t>0$, there is a constant $c_1=c_1(t)>1$ such that
 \beqnn
 \ar\ar
\sup_{s\in(0,t],\,x>c_1}[\ln x]^{-1}[x^{-1}\gamma_0(s,x)-H(s,x)]
<\infty,
 \eeqnn
then $\mbf{P}\{\tau_\infty<\infty\}=0$.
\end{itemize}
\etheorem
The conditions in Theorem \ref{t1.3} generalize those in \cite[Theorems 2.3 (i) and 2.8 (i)]{LYZh} for the case when $\gamma_0(s,x)=\gamma_0(x)$ for some Borel function on $[0,\infty)$ and $b_1(s)=b_2(s)=1$ for all $s>0$.

Throughout the remainder of this section, we assume that $\int_0^1 z\mu(\dd z)=\infty$.
We now proceed to study the extinction and explosion
conditions for the solution to \eqref{ii1.5}.

\btheorem\label{t1.4}
\begin{itemize}
\item[{\normalfont(i)}]
If for some $t_0>0$, there are constants $\delta>1$ and $0<c_0,\delta_0<1$ and a strictly positive function $d$ on $(0,t_0]$ such that
 \beqnn
\inf_{s\in(0,t_0],\,0<x<c_0}[\ln x^{-1}]^{-\delta}d(s)^{-1}[-x^{-1}\gamma_0(s,x)+\delta_0H(s,x)]
>0,
 \eeqnn
then
$\mbf{P}\{\tau_0<\infty\}>0$
for all small enough $X_0>0$.
\item[{\normalfont(ii)}]
If for some $\tilde{t}_0>0$, there are constants $\tilde{\delta},\tilde{c}_0,\tilde{\delta}_0>1$ and a strictly positive function $\tilde{d}$ on $(0,\tilde{t}_0]$ such that
 \beqnn
\inf_{s\in(0,\tilde{t}_0],\,x>\tilde{c}_0}[\ln x ]^{-\tilde{\delta}}\tilde{d}(s)^{-1}[x^{-1}\gamma_0(s,x)- \tilde{\delta}_0H(s,x)]
>0,
 \eeqnn
then
$\mbf{P}\{\tau_\infty<\infty\}>0$ for all large enough $X_0>0$.
\end{itemize}
\etheorem
The conditions in Theorem \ref{t1.4} reduce to those in \cite[Theorems 2.3(i) and 2.8(i)]{LYZh} when $\gamma_0(s,x) = \gamma_0(x)$ for some Borel function on $[0,\infty)$ and $b_1(s) = b_2(s) = 1$ for all $s > 0$.

Note that the parameter ranges in Theorems \ref{t1.3} and \ref{t1.4} are mutually exclusive. Regarding extinction behavior, Theorem \ref{t1.3}   (i) establishes non-extinction under a boundedness condition near zero, while Theorem \ref{t1.4}   (i) yields positive extinction probability under a strict positivity condition with parameters $\delta > 1$. Similarly, for explosion behavior, Theorem \ref{t1.3}   (ii) ensures non-explosion through a boundedness criterion near infinity, whereas Theorem \ref{t1.4}   (ii) establishes positive explosion probability under a strict positivity requirement with $\tilde{\delta} > 1$.

This opposition in parameter ranges, together with the complementary probability conclusions of almost sure non-extinction versus possible extinction and almost sure non-explosion versus possible explosion, suggests that our conditions may provide effective threshold criteria for characterizing the qualitative behavior of the process.

We now establish the extinction and explosion conditions that hold for all $X_0>0$.

\btheorem\label{t1.4b}
\begin{itemize}
\item[{\normalfont(i)}]
Suppose that the assumption in Theorem \ref{t1.4} (i) holds and
there is a constant $b>0$ so that
one of the following holds for all $0<a<X_0$:
\begin{itemize}
\item[{\normalfont(ia)}]
$\inf\limits_{x\in[a,X_0+b]}\gamma_1(x)>0$,
$b_1(s)>0$ for all $0<s\le t_0$, $\sup\limits_{x\in[a,X_0+b],s\in(0,t_0]} \gamma_0(s,x)/b_1(s)<\infty$;
\item[{\normalfont(ib)}]
$\inf\limits_{x\in[a,X_0+b]}\gamma_2(x)>0$,
$b_2(s)>0$ for all $0<s\le t_0$, $\sup\limits_{x\in[a,X_0+b],s\in(0,t_0]} \gamma_0(s,x)/b_2(s)<\infty$.
\end{itemize}
Then
$\mbf{P}\{\tau_0<\infty\}>0$  for all  $X_0>0$.
\item[{\normalfont(ii)}]
Suppose that the assumption in Theorem \ref{t1.4} (ii) holds and
there is a constant $0<a<X_0$ so that
one of (ia)--(ib) holds for all $b>0$.
Then
$\mbf{P}\{\tau_\infty<\infty\}>0$ for all  $X_0>0$.
\end{itemize}
\etheorem
When $\gamma_0(s,x)=\gamma_0(x)$ for some Borel function on $[0,\infty)$ and $b_1(s)=b_2(s)=1$ for all $s>0$,
Theorem \ref{t1.4b} (i) and (ii) correspond to
\cite[Corollary 2.5 and Proposition 2.6 (i)]{LYZh}
and \cite[Remark 2.10 and Proposition 2.11]{LYZh}, respectively.

We now provide conditions for almost sure extinction.
\btheorem\label{t1.4c}
Let
 \beqnn
 H_\rho(s,x)
  \ar=\ar
-x^{\rho-1}\gamma_0(s,x)+(1-\rho)b_1(s)x^{\rho-2}\gamma_1(x) \cr
 \ar\ar
+(1-\rho)b_2(s)\gamma_2(x)
 \int_0^\infty z^2\mu(\dd z)
\int_0^1(x+zu)^{\rho-2}(1-u)\dd u.
 \eeqnn
Then
$\mbf{P}\{\tau_0<\infty\}=1$ if there are a constant $0<\rho<1$
and a strictly positive function $d$
with $d((0,\infty))=\infty$ such that
$\inf_{s>0,\,0<x<b}H_\rho(s,x)/d(s)>0$
for all $b>0$.
\etheorem

Finally, we study the first passage probabilities. These results also generalize those given in
\cite[Propositions 2.6 (i) and 2.11]{LYZh}, which correspond to the special case where $\gamma_0(s,x)=\gamma_0(x)$ for some Borel function on $[0,\infty)$ and $b_1(s)=b_2(s)=1$ for all $s>0$.

\bproposition\label{t1.5}
Suppose that $0<a<X_0<b$.
Then
$\mbf{P}\{\tau_a^-<\infty\}>0$ and $\mbf{P}\{\tau_b^+<\infty\}>0$
if there is a constant $t_0>0$ such that one of the following holds:
\begin{itemize}
\item[{\normalfont(i)}]
$\inf_{x\in[a,b]}\gamma_1(x)>0$,
$b_1(s)>0$ for all $0<s\le t_0$
and $\sup_{x\in[a,b],s\in(0,t_0]} |\gamma_0(s,x)|/b_1(s)<\infty$;
\item[{\normalfont(ii)}]
$\inf_{x\in[a,b]}\gamma_2(x)>0$,
$b_2(s)>0$ for all $0<s\le t_0$
and $\sup_{x\in[a,b],s\in(0,t_0]}  |\gamma_0(s,x)|/b_2(s)<\infty$.
\end{itemize}
\eproposition

\subsection{Proofs}

We begin with a lemma that will be used repeatedly.

\blemma\label{t2.1}
Suppose that \eqref{3.0} holds.
Let $b>a>0$ and function $g\in C^2([a,\infty))$
satisfy $\sup_{x\ge a}|g'(x)|\wedge \sup_{x\ge a}|g(x)|<\infty$.
Then for the stopping time $\varsigma$ satisfying $a\le X_s\le b$  for     $s<\varsigma$, we have
 \beqnn
\mbf{E}[g(X_{t\wedge\varsigma})]
=g(X_0)+\int_0^t\mbf{E}\Big[\mathcal{L}_sg(X_s)1_{\{s\le \varsigma\}}\Big]\dd s.
 \eeqnn
\elemma
\proof
By It\^o's formula and \eqref{ii1.5} we have
 \beqlb\label{3.10}
g(X_{t\wedge\varsigma})
 \ar=\ar
g(X_0)+\int_0^{t\wedge\varsigma}\mathcal{L}_sg(X_s)\dd s
+M_1(t\wedge\varsigma)+M_2(t\wedge\varsigma),
 \eeqlb
where
 \beqnn
M_1(t):=\int_0^t\int_0^{\gamma_1(X_{s-})}g'(X_{s-})W(\dd s,\dd u)
 \eeqnn
and
 \beqnn
M_2(t):=\int_0^t\int_0^\infty\int_0^{\gamma_2(X_{s-})}
\Delta_zg(X_{s-})\tilde{N}(\dd s,\dd z,\dd u),
 \eeqnn
here, the operator $\Delta_z$ is defined at the end of Section 1.2.  Under \eqref{3.0} and the conditions of the lemma,
 \beqnn
\gamma_1(u)|g'(u)|^2\le
\sup_{u\in[a,b]}[\gamma_1(u)|g'(u)|^2]=:c_1
 \eeqnn
and then
 \beqnn
\int_0^t\int_0^{\gamma_1(X_{s-})}|g'(X_{s-})|^2b_1(s)\dd s\dd u\le c_1
b_1((0,t])<\infty.
 \eeqnn
Thus, $t\mapsto M_1(t\wedge\varsigma)$ is a martingale.
Similarly,
 \beqnn
\gamma_2(u)|\Delta_zg(u)|^2\le
z^2\sup_{u\in[a,b]}\gamma_2(u)
\sup_{a\le v\le b+1}|g'(v)|^2=:z^2c_2,\quad 0<z\le1
 \eeqnn
and
 \beqnn
\gamma_2(u)|\Delta_zg(u)|\le
2z\sup_{u\in[a,b]}\gamma_2(u)[\sup_{v\ge a}|g'(v)|\wedge
\sup_{v\ge a}|g(v)|]=:zc_3,\quad z>1.
 \eeqnn
Thus
 \beqnn
\int_0^{t\wedge\varsigma}\int_0^1\int_0^{\gamma_2(X_{s-})} |\Delta_zg(u)|^2b_2(s)\dd s\mu(\dd z)\dd u\le c_2b_2((0,t])\int_0^1z^2\mu(\dd z)
<\infty
 \eeqnn
and
 \beqnn
\int_0^{t\wedge\varsigma}\int_1^\infty\int_0^{\gamma_2(X_{s-})} |\Delta_zg(u)|b_2(s)\dd s\mu(\dd z)\dd u\le c_3b_2((0,t])\int_1^\infty z\mu(\dd z)<\infty.
 \eeqnn
Thus $t\mapsto M_2(t\wedge\varsigma)$ is also a martingale.
Taking expectations on both sides of \eqref{3.10} yields the assertion.
\qed

The proof of Theorem \ref{t1.3} adapts the approach of \cite[Theorem 2.25]{chen04} for establishing the uniqueness of the $q$-process to our setting.
The key test functions $g$ are of logarithmic type, chosen to diverge to infinity at zero (for non-extinction) and at infinity (for non-explosion).
The crucial step is to show that $\mathcal{L}_sg$ is bounded above by $g$ (see \eqref{3.11} and \eqref{3.12}).

\noindent{\it Proof of Theorem \ref{t1.3}.}
(i)
For $n\ge1$ let nonnegative function $g_n\in C^2((0,\infty))$
satisfy $g_n(x):=1+\ln n+\ln x^{-1}$ for $0<x\le n$, $g_n''(x)=0$
for $x\ge n+1$, and $g_n''(x)\le x^{-2}$ for all $x>0$.
Then
 \beqnn
\int_0^\infty z^2\mu(\dd z)\int_0^1g_n''(x+zu)(1-u)\dd u
\le\int_0^\infty z^2\mu(\dd z)\int_0^1(x+zu)^{-2}(1-u)\dd u,\quad
0<x<n.
 \eeqnn
It follows from \eqref{3.0a} and
\eqref{3.aaa} that there is a constant $C_{n,T}>0$ so that
 \beqlb\label{3.11}
\mathcal{L}_sg_n(x)
\le
-x^{-1}\gamma_0(s,x)+
H(s,x)
\le
C_{n,T}g_n(x),\qquad s\in(0,T],~x\in(0,n].
 \eeqlb
By Lemma \ref{t2.1},
 \beqnn
\mbf{E}[g_n(X_{t\wedge\tau_{m,n}})]\le
g_n(X_0)+C_{n,T}\int_0^t\mbf{E}[g_n(X_{s\wedge\tau_{m,n}})]\dd s
 \eeqnn
for all $t\in(0,T]$.
Using Gronwall's inequality we obtain
 \beqnn
\mbf{E}[g_n(X_{t\wedge\tau_{m,n}})]
\le
g_n(X_0) \e^{C_{n,T} t},\qquad t\in(0,T].
 \eeqnn
Taking $m\to\infty$ and using Fatou's lemma we have
 \beqnn
\mbf{E}[g_n(X_{t\wedge\tau_0\wedge\tau_n^+})]
\le
\mbf{E}[\lim_{m\to\infty}g_n(X_{t\wedge\tau_{m,n}})]
\le\liminf_{m\to\infty}\mbf{E}[g_n(X_{t\wedge\tau_{m,n}})]
\le g_n(X_0) \e^{C_{n,T} t},~~ t\in(0,T].
 \eeqnn
Since $\lim_{x\to0}g_n(x)=0$,
it follows that $\mbf{P}\{\tau_0\le t\wedge\tau_n^+\}=0$.
Taking $t,n\to\infty$ one has $\tau_0>\tau_\infty$ almost surely.
The conclusion then follows from Definition \ref{Def of weak solu}.

(ii)
For $m\ge1$ and nonnegative function $g_m\in C^2((0,\infty))$
satisfying $g_m(x):=1+\ln m+\ln x$ for $x\ge m^{-1}$, $g_n''(x)=0$
for $0<x<(m+1)^{-1}$, and $g_n''(x)\ge -x^{-2}$ for all $x>0$.
Following the same argument as in (i), for all $T>0$ and $m>1$, there
is a constant $\tilde{C}_{m,T}>0$ such that
 \beqlb\label{3.12}
\mathcal{L}_sg_m(x)
\le
x^{-1}\gamma_0(s,x)-H(s,x)
\le
\tilde{C}_{m,T}g_m(x),\qquad s\in(0,T],~x\ge m^{-1}.
 \eeqlb
and then
$\mbf{E}[g_m(X_{t\wedge\tau_{m,n}})]
\le
g_m(X_0) \e^{\tilde{C}_{m,T} t}$ for all $t\in(0,T]$.
Taking $n\to\infty$ we get
 \beqnn
\mbf{E}[g_m(X_{t\wedge\tau_\infty\wedge\tau_{1/m}^-})]
\le g_m(X_0) \e^{\tilde{C}_{m,T} t},~~ t\in(0,T].
 \eeqnn
Since $\lim_{x\to\infty}g_m(x)=0$, we have
$\mbf{P}\{\tau_\infty\le t\wedge\tau_{1/m}^-\}=0$.
Then letting $t,m\to\infty$,  we obtain the desired result.
\qed

The proof of Theorem \ref{t1.4} follows the strategy of  \cite[Theorem 2.27]{chen04} for demonstrating the non-uniqueness of the $q$-process.
In our setting, we introduce a new family of test functions $g$, defined in Lemma \ref{t1.7}, which are compositions of logarithmic and exponential functions. These functions are designed to remain bounded. The  main challenge is to show that $\mathcal{L}_sg$ is bounded below by $g$ near zero and infinity, respectively; see \eqref{2.12} and \eqref{2.12b}.

We begin with the following estimate
on $\mathcal{L}_sg$.

\blemma\label{t1.7}
\begin{itemize}
\item[{\normalfont(i)}]
Given $\lambda>1$ and $\rho>0$, let $g(x)=\e^{-\lambda h(x)}$
where $h\in C^2((0,\infty))$ is a nonnegative and non-decreasing function with
$h(x)=(\ln x^{-1})^{-\rho}$ for  $0<x<1/2$.
Then there is a constant $c\in(0,1/2)$ such that
 \beqnn
\mathcal{L}_sg(x)
\ge
\lambda \rho g(x)(\ln x^{-1})^{-\rho-1}
\big[-x^{-1}\gamma_0(s,x)+\delta_0H(s,x)\big],\qquad
0<x<c,~s>0,
 \eeqnn
where $0<\delta_0<1$ is the constant given in Theorem \ref{t1.4}(i).
\item[{\normalfont(ii)}]
Given $\lambda>1$ and $\rho>0$, let $g(x)=\e^{-\lambda h(x)}$,
where $h\in C^2((0,\infty))$ is a nonnegative and non-increasing function with
$h(x)=(\ln x)^{-\rho}$ for all $x>2$.
Then there is a constant $c>2$ such that
 \beqnn
\mathcal{L}_sg(x)
\ge
\lambda \rho g(x)(\ln x^{-1})^{-\rho-1}
\big[x^{-1}\gamma_0(s,x)-\tilde{\delta}_0H(s,x)\big],\qquad
x>c,~s>0,
 \eeqnn
where $\tilde{\delta}_0>1$ is the constant given in Theorem \ref{t1.4}(ii).
\end{itemize}
\elemma
\proof
Since the proofs are similar, we prove only part  (i).
For $0<x<3^{-2}$,
 \beqlb\label{3.5}
g'(x)[\lambda\rho g(x)]^{-1}=-(\ln x^{-1})^{-\rho-1} x^{-1}
 \eeqlb
and
 \beqlb\label{3.6}
g''(x)[\lambda\rho g(x)]^{-1}
=
(\ln x^{-1})^{-\rho-1} x^{-2}
-(\rho+1)(\ln x^{-1})^{-\rho-2} x^{-2}
+\lambda \rho (\ln x^{-1})^{-2\rho-2} x^{-2}.
 \eeqlb
Moreover, there are constants $0<\delta_1<1-\delta_0$ and $0<c_1<3^{-3}$ so that for $
0<x<c_1$,
 \beqlb\label{3.8}
g''(x)[\lambda\rho g(x)]^{-1}\ge(1-\delta_1)(\ln x^{-1})^{-\rho-1} x^{-2},
 \eeqlb
and for $0<x<2c_1$,
 \beqnn
h''(x)=-\rho(\ln x^{-1})^{-\rho-1} x^{-2}[1-(\rho+1)(\ln x^{-1})^{-1}]
\le -\rho(1-\delta_1/2)(\ln x^{-1})^{-\rho-1} x^{-2}.
 \eeqnn
Since $\e^{x}-1\ge x$ for all $x\in\mbb{R}$,   applying Taylor expansion twice we have
 \beqnn
 \ar\ar
z^2\int_0^1 g''(x+zu)(1-u)\dd u
=D_zg(x) \cr
 \ar\ge\ar
-\lambda g(x)D_zh(x)
=
-\lambda g(x)z^2\int_0^1h''(x+zu)(1-u)\dd u \cr
 \ar\ge\ar
(1-\delta_1/2)\rho\lambda g(x)z^2\int_0^1(\ln (x+zu)^{-1})^{-\rho-1} (x+zu)^{-2}(1-u)\dd u \cr
 \ar\ge\ar
(1-\delta_1/2)\rho\lambda g(x)(\ln x^{-1})^{-\rho-1} G(x,z),\qquad 0<x,z<c_1,
 \eeqnn
where
 \beqnn
G(x,z):=z^2\int_0^1(x+zu)^{-2}(1-u)\dd u
=
-\ln(1+x^{-1}z)+x^{-1}z
 \eeqnn
by Taylor expansion again.
By the assumptions on $g$,
 \beqnn
\int_{c_1}^\infty D_zg(x)\mu(\dd z)
\ge
-g(x)\int_{c_1}^\infty\mu(\dd z),\qquad 0<x<c_1.
 \eeqnn
It follows that for $0<x<c_1$,
\beqlb\label{3.7}
\ar\ar\int_0^\infty z^2\mu(\dd z)\int_0^1 g''(x+zu)(1-u)\dd u \cr
 \ar\ge\ar
(1-\delta_1/2)\rho\lambda(\ln x^{-1})^{-\rho-1}g(x)
\int_0^{c_1}G(x,z)\mu(\dd z)
-g(x)\int_{c_1}^\infty\mu(\dd z).
 \eeqlb
By Fatou's lemma,
 \beqnn
\lim_{x\to0}x\int_0^{c_1}G(x,z)\mu(\dd z)
\ge
\int_0^{c_1}\lim_{x\to0}xG(x,z)\mu(\dd z)
=\int_0^{c_1}z\mu(\dd z)=\infty
 \eeqnn
and
 \beqnn
\lim_{x\to0}(\ln x^{-1})^{-\rho-1}
\int_0^{c_1}G(x,z)\mu(\dd z)
\ge
\int_0^{c_1}\lim_{x\to0}(\ln x^{-1})^{-\rho-1}G(x,z)\mu(\dd z)=\infty.
 \eeqnn
Then there is a constant $0<c_2<c_1$ so that
 \beqlb\label{3.9}
2^{-1}\delta_1\rho(\ln x^{-1})^{-\rho-1}\int_0^{c_1}G(x,z)\mu(\dd z)
\ge
\int_{c_1}^\infty\mu(\dd z),\quad
0<x<c_2,
 \eeqlb
and
 \beqnn
(1-\delta_0-\delta_1)x\int_0^{c_1}G(x,z)\mu(\dd z)
\ge x\int_{c_1}^\infty x^{-1}z\mu(\dd z)
\ge x\int_{c_1}^\infty G(x,z)\mu(\dd z),\quad
0<x<c_2,
 \eeqnn
which implies
 \beqnn
\delta_0\int_0^\infty G(x,z)\mu(\dd z)
\le
(1-\delta_1)\int_0^{c_1} G(x,z)\mu(\dd z) ,\quad
0<x<c_2.
 \eeqnn
Combining this with \eqref{3.7} and \eqref{3.9} one gets
 \beqnn
 \ar\ar
\int_0^\infty z^2\mu(\dd z)\int_0^1 g''(x+zu)(1-u)\dd u \cr
 \ar\ge\ar
(1-\delta_1)\rho\lambda(\ln x^{-1})^{-\rho-1}g(x)
\int_0^{c_1}G(x,z)\mu(\dd z) \cr
 \ar\ge\ar
\delta_0\rho\lambda(\ln x^{-1})^{-\rho-1}g(x) \int_0^\infty G(x,z)\mu(\dd z)
 \eeqnn
for all $0<x<c_2$.
Combining this with \eqref{3.5} and \eqref{3.8} yields
the assertion with $c:=c_2$.
\qed

Now we are ready to prove Theorem \ref{t1.4}.

\noindent{\it Proof of Theorem \ref{t1.4}.}
(i)
Take $0<\rho<\delta-1$.
Let $g$ be the function defined in Lemma \ref{t1.7} (i).
Let $d_0:=d((0,t_0])$ and
$d(s)=d(t_0)$ for all $s\ge t_0$.
Then $d_0>0$ and $d((0,\infty))=\infty$.
By Lemma \ref{t1.7}   (i) and the assumption, there is a constant $C_0>0$ such that
 \beqlb\label{3.3}
\mathcal{L}_sg(x)\ge \lambda g(x)\rho C_0 (\ln x^{-1})^{\delta-\rho-1}d(s)
,\qquad 0<x<c_0\wedge c,\,0<s<t_0.
 \eeqlb
Observe that there are constants $0<c_1<c_0\wedge c$ and $\bar{X}_0>0$ such that
 \beqlb\label{2.10a}
(\ln c_1^{-1})^{-\rho}\ge 2(\ln \bar{X}_0^{-1})^{-\rho},\quad
d_0\rho C_0 (\ln c_1^{-1})^{\delta-\rho-1}\ge 2(\ln \bar{X}_0^{-1})^{-\rho}.
 \eeqlb
In the following, we take $0<X_0<\bar{X}_0$ and
let $\delta_1:=2(\ln X_0^{-1})^{-\rho}/d_0$.
Then by \eqref{2.10a},
 \beqlb\label{2.10}
(\ln c_1^{-1})^{-\rho}\ge \delta_1d_0,\qquad
\rho C_0 (\ln c_1^{-1})^{\delta-\rho-1}\ge \delta_1.
 \eeqlb
Choose $\lambda>0$ sufficiently large such that
 \beqlb\label{2.11}
1-2\e^{-\lambda \delta_1d_0/2}>0.
 \eeqlb
By \eqref{3.3}, for all $0<x<c_1$ and $0<s\le t_0$,
 \beqlb\label{2.12}
\mathcal{L}_sg(x)\ge \lambda g(x)\rho C_0 (\ln x^{-1})^{\delta-\rho-1}d(s)
\ge
\lambda g(x)\rho C_0 (\ln c_1^{-1})^{\delta-\rho-1}d(s)
\ge
\delta_1\lambda g(x) d(s).
 \eeqlb

Let $\gamma_{m}:=\tau_{1/m}^-\wedge\tau_{c_1}^+\wedge t_0$.
By Lemma \ref{t2.1},
 \beqnn
\mbf{E}[g(X_{t\wedge\gamma_m})]
=g(X_0)
+\int_0^t\mbf{E}[\mathcal{L}_sg(X_s)1_{\{s\le\gamma_{m}\}}]\dd s.
 \eeqnn
Using integration by parts and \eqref{2.12}, we have
 \beqnn
 \ar\ar
\e^{-\lambda\delta_1d((0,t])}\mbf{E}[g(X_{t\wedge\gamma_{m}})] \cr
 \ar=\ar
g(X_0)-\lambda\delta_1\int_0^td(s)\e^{-\lambda\delta_1d((0,s])}\mbf{E}[g(X_{s\wedge\gamma_{m}})]\dd s
+\int_0^t\e^{-\lambda\delta_1d((0,s])}\mbf{E}[\mathcal{L}_sg(X_s)1_{\{s\le\gamma_{m}\}}]\dd s \cr
 \ar\ge\ar
g(X_0)+\mbf{E}\Big[\lambda\delta_1g(X_{\gamma_{m}})
\int_{\gamma_{m}}^{t\vee\gamma_{m}}\e^{-\lambda\delta_1d((0,s])}d(s)\dd s \Big].
 \eeqnn
Letting $t\to\infty$ and using the fact $d((0,\infty))=\infty$,  we obtain
 \beqlb\label{2.14}
g(X_0)
\le
\mbf{E}\Big[g(X_{\gamma_{m}})
\e^{-\lambda\delta_1d((0,\gamma_{m}])}\Big].
 \eeqlb
Letting $m\to\infty$ and applying the dominated convergence theorem, we obtain
 \beqlb\label{2.4}
\begin{array}{ll}
g(X_0)
 \le
\mbf{E}\Big[g(X_{\tau_0\wedge\tau_{c_1}^+\wedge t_0})
\e^{-\lambda \delta_1d((0,\tau_0\wedge\tau_{c_1}^+\wedge t_0])}\Big] \\
\qquad~~\,\le
\mbf{E}\Big[g(X_{\tau_0\wedge\tau_{c_1}^+\wedge t_0})\e^{-\lambda \delta_1d((0,\tau_0\wedge\tau_{c_1}^+\wedge t_0])}\Big(1_{\{\tau_0\le\tau_{c_1}^+\wedge t_0\}}+1_{\{\tau_{c_1}^+\le\tau_0\wedge t_0\}}
+1_{\{t_0\le\tau_0\wedge \tau_{c_1}^+\}}\Big)\Big] \\
\qquad~~\,\le
\mbf{P}\{\tau_0\le t_0\}+g(c_1)+
\e^{-\lambda \delta_1 d_0}.
\end{array}
 \eeqlb
By \eqref{2.10} we have $g(c_1)\le g(X_0)^2$
and $\e^{-\lambda \delta_1 d_0}=g(X_0)^2$,
and then using \eqref{2.11} we obtain
 \beqnn
\mbf{P}\{\tau_0\le t_0\}
\ge g(X_0)-2g(X_0)^2
\ge g(X_0)[1-2g(X_0)]>0.
 \eeqnn
This proves part (i).

(ii)
Take $0<\tilde{\rho}<\tilde{\delta}-1$.
Let $\tilde{d}_0:=\tilde{d}((0,\tilde{t}_0])$ and
$\tilde{d}(s)=\tilde{d}(\tilde{t}_0)$ for all $s\ge \tilde{t}_0$.
Let $g$ be the function defined in Lemma \ref{t1.7} (ii).
Then by Lemma \ref{t1.7} (ii) and the assumption, similar to \eqref{2.12},
there are large enough constants $\tilde{c}_1>\tilde{c}_0$ and $\bar{X}_0>0$ such that for  $X_0>\bar{X}_0$,
 \beqlb\label{2.13}
(\ln \tilde{c}_1)^{-\tilde{\rho}}\ge 2(\ln \bar{X}_0)^{-\tilde{\rho}}
\ge 2(\ln X_0)^{-\tilde{\rho}}=:\tilde{\delta}_1\tilde{d}_0
 \eeqlb
and
 \beqlb\label{2.12b}
\mathcal{L}_sg(x)\ge
\lambda \tilde{\delta}_1g(x) \tilde{d}(s),\qquad x>\tilde{c}_1,~0<s\le \tilde{t}_0.
 \eeqlb
Take $\lambda>0$ sufficiently large such that
$1-2\e^{-\lambda (\ln X_0)^{-\tilde{\rho}}}>0$.
By Lemma \ref{t2.1}
and following an argument similar to the derivation of \eqref{2.4}, we have
 \beqnn
g(X_0)
 \ar\le\ar
\mbf{E}\Big[g(X_{\tau_\infty\wedge\tau_{\tilde{c}_1}^-\wedge \tilde{t}_0})
\e^{-\lambda \tilde{\delta}_1\tilde{d}((0,\tau_\infty\wedge\tau_{\tilde{c}_1}^-\wedge \tilde{t}_0])}\Big] \cr
 \ar\le\ar
\mbf{E}\Big[g(X_{\tau_\infty\wedge\tau_{\tilde{c}_1}^-\wedge \tilde{t}_0})\e^{-\lambda\tilde{\delta}_1\tilde{d}((0,\tau_\infty\wedge\tau_{\tilde{c}_1}^-\wedge \tilde{t}_0])}\Big(1_{\{\tau_\infty\le\tau_{\tilde{c}_1}^-\wedge \tilde{t}_0\}}+1_{\{\tau_{\tilde{c}_1}^-\le\tau_\infty\wedge \tilde{t}_0\}}
+1_{\{\tilde{t}_0\le\tau_\infty\wedge \tau_{\tilde{c}_1}^-\}}\Big)\Big] \cr
 \ar\le\ar
\mbf{P}\{\tau_\infty\le \tilde{t}_0\}+g(\tilde{c}_1)+
\e^{-\lambda\tilde{\delta}_1 \tilde{d}_0}.
 \eeqnn
It follows from \eqref{2.13} that
 \beqnn
\mbf{P}\{\tau_\infty\le \tilde{t}_0\}
\ge g(X_0)[1-2g(X_0)]>0,
 \eeqnn
which establishes the desired conclusion.
\qed

For the proof of Theorem \ref{t1.4b} we employ the same
test function as in  Theorem \ref{t1.4}.
The key step is to verify that $\mathcal{L}_s g$ is bounded below by $g$ on an interval containing $(0, X_0)$ for the extinction behavior (see \eqref{3.13aa} and \eqref{3.13ab}) and on $(X_0, \infty)$ for the explosion behavior.

\noindent{\it Proof of Theorem \ref{t1.4b}.}
(i)
Let  $g$
be the function defined in Lemma \ref{t1.7} (i).
Since the proofs for  Conditions (ia) and (ib) are similar, we present only the proof under condition (ia).
Recall the constant $t_0$ from Theorem \ref{t1.4}(i).
For simplicity, we assume that $X_0+b<1$.
Define $\bar{d}(s):=d(s)\wedge b_1(s)$.
Let $d_0:=\bar{d}((0,t_0])$ and
$\bar{d}(s)=\bar{d}(t_0)$ for all $s\ge t_0$.
By \eqref{3.3}, there exists a constant $c_1>0$ so that
 \beqlb\label{3.13aa}
\mathcal{L}_sg(x)\ge 2\lambda d_0^{-1} g(x) d(s)
,\qquad 0<x<c_1,\,0<s<t_0.
 \eeqlb
From \eqref{3.5} and \eqref{3.6}, under  assumption (ia),
there exists  $\lambda_1>0$
such that for all $\lambda>\lambda_1$, $
c_1\le x\le X_0+b$ and $0<s<t_0$,
 \beqlb\label{3.13ab}
\mathcal{L}_sg(x)
\ar\ge\ar
b_1(s)\lambda \rho g(x)(\ln x^{-1})^{-\rho-1}
\big[-x^{-1}\gamma_0(s,x)/b_1(s) \cr
\ar\ar\,\,\,\qquad
+x^{-2}\gamma_1(x)
(1-(\rho+1)(\ln x^{-1})^{-1}+\lambda \rho (\ln x^{-1})^{-\rho-1}) \big] \cr
\ar\ge\ar
2\lambda d_0^{-1}  g(x) b_1(s).
 \eeqlb
Following an argument similar to the derivation of \eqref{2.4}, we obtain for all $\lambda > \lambda_1$,
 \beqnn
g(X_0)
 \ar\le\ar
\mbf{E}\Big[g(X_{\tau_0\wedge\tau_{X_0+b}^+\wedge t_0})
\e^{-2\lambda d_0^{-1} \bar{d}((0,\tau_0\wedge\tau_{X_0+b}^+\wedge t_0])}\Big] \cr
 \ar\le\ar
\mbf{P}\{\tau_0\le t_0\}+g(X_0+b)+
\e^{-2\lambda(\ln X_0^{-1})^{-\rho} } \cr
 \ar\le\ar
\mbf{P}\{\tau_0\le t_0\}+g(X_0+b)+g(X_0)^2,
 \eeqnn
 where the second inequality follows by analyzing the value of the stopped process $X_t$ at the terminal time and bounding the exponential term.
Now, we choose $\lambda_2 > \lambda_1$ such that for all $\lambda > \lambda_2$,
 \beqnn
1-\e^{-\lambda [(\ln (X_0+b)^{-1})^{-\rho}-(\ln X_0^{-1})^{-\rho}]}
-\e^{-\lambda (\ln X_0^{-1})^{-\rho}}>0.
 \eeqnn
This condition ensures that the combination of terms involving $g(X_0+b)$ and $g(X_0)^2$ does not dominate $g(X_0)$.
It then follows that
 \beqnn
\mbf{P}\{\tau_0\le t_0\}\ge g(X_0)-g(X_0+b)-g(X_0)^2>0.
 \eeqnn
This proves assertion (i).

(ii) The proof of part (ii) is a modification of that of part (i) and is omitted here.
\qed

The proof of Theorem \ref{t1.4c} employs the same technique of using a composition of a power function and an exponential function as the key test function. The crucial estimate for the operator $\mathcal{L}_s$ is given by \eqref{3.13b}.

\noindent{\it Proof of Theorem \ref{t1.4c}.}
Let $g(x)=\e^{-\lambda x^{\rho}}$ for $x,\lambda>0$.
Under the stated assumptions, we have for some $c_b > 0$ and all $b > 0$,
 \beqlb\label{3.13b}
\mathcal{L}_sg(x)
\ge \lambda\rho H_{\rho}(s,x)
\ge \lambda\rho c_b d(s),\qquad 0<x\le b,~s>0
 \eeqlb
Following an argument similar to that for \eqref{2.4}, we obtain
 \beqnn
g(X_0)
 \ar\le\ar
\mbf{E}\Big[g(X_{\tau_0\wedge\tau_{b}^+})
\e^{-\lambda\rho c_b d((0,\tau_0\wedge\tau_{b}^+])}1_{\{\tau_0\wedge\tau_{b}^+<\infty\}}\Big] \cr
 \ar\le\ar
\mbf{E}\Big[g(X_{\tau_0\wedge\tau_{b}^+})
\e^{-\lambda\rho c_b d((0,\tau_0\wedge\tau_{b}^+])}\Big(1_{\{\tau_0<\tau_{b}^+\}}
+1_{\{\tau_{b}^+<\tau_0\}}\Big)
\Big] \cr
 \ar\le\ar
\mbf{P}\{\tau_0<\tau_{b}^+\}+g(b)
\le
\mbf{P}\{\tau_0<\infty\}+g(b).
 \eeqnn
Letting $b\to\infty$ and then $\lambda\to 0$ completes the proof.
\qed

The proof of Proposition \ref{t1.5} follows a similar method, with exponential functions serving as the key test functions. The required estimates for the operator $\mathcal{L}_s$ are provided in \eqref{3.13ba} and \eqref{3.13bb}.

\noindent{\it Proof of Proposition \ref{t1.5}.}
Since the proofs under conditions (i) and (ii) are similar, we only present the proof for condition (i).
Let
$g(x)=\e^{-\lambda x}$ for $\lambda,x>0$.
Let $d_0:=b_1((0,t_0])$, and $d(s)=b_1(s)$ for $0<s\le t_0$ and $d(s)=b_1(t_0)$ for $s>t_0$.
Then there is a constant $\lambda_1>0$ such that
 \beqnn
2^{-1}\lambda_1 \inf_{x\in[a,b]}\gamma_1(x)
\ge
\sup_{a\le x\le b,s\in(0,t_0]}|\gamma_0(s,x)| /b_1(s),
\quad\sqrt{\lambda_1}\inf_{x\in[a,b]}\gamma_1(x)\ge2,
 \eeqnn
and for $\lambda=\lambda_1$,
 \beqlb\label{3.4}
g(X_0)-g(b)-\e^{-\lambda^{3/2}d_0}>0.
 \eeqlb
Consequently, for $\lambda = \lambda_1$ and all $a \le x \le b$, $0 < s \le t_0$,
 \beqlb\label{3.13ba}
\mathcal{L}_sg(x)
\ar = \ar
\lambda\e^{-\lambda x}
\Big[-\gamma_0(s,x)+\lambda \gamma_1(x)b_1(s)
+\gamma_2(x)b_2(s)\int_0^\infty [\lambda^{-1}\e^{-\lambda z}-\lambda^{-1}+z]\mu(\dd z)\Big] \cr
\ar\ge\ar
\lambda\e^{-\lambda x}b_1(s)[-|\gamma_0(s,x)|/b_1(s)+2^{-1}\lambda \gamma_1(x)]
+2^{-1}\lambda^2 \gamma_1(x)b_1(s)g(x) \cr
\ar\ge \ar
\lambda^{3/2}d(s)g(x).
 \eeqlb
Let $\gamma_{a,b}=\tau_a^-\wedge\tau_b^+\wedge t_0$.
Then, for $\lambda = \lambda_1$, an argument similar to that for \eqref{2.4} yields
 \beqnn
g(X_0)
 \ar\le\ar
\mbf{E}\Big[g(X_{\gamma_{a,b}})
\e^{-\lambda^{3/2}d((0,\gamma_{a,b}])}
\Big] \cr
 \ar\le\ar
\mbf{E}\Big[g(X_{\gamma_{a,b}})
\e^{-\lambda^{3/2}d((0,\gamma_{a,b}])}
\Big(1_{\{\tau_a^-\le\tau_b^+\wedge t_0\}}
+1_{\{\tau_b^+\le\tau_a^-\wedge t_0\}}
+1_{\{t_0\le\tau_a^-\wedge \tau_b^+\}}\Big)\Big] \cr
 \ar\le\ar
g(a)\mbf{P}\{\tau_a^-\le\tau_b^+\wedge t_0\}
+g(b)+
\e^{-\lambda^{3/2}d_0}.
 \eeqnn
Using \eqref{3.4}, we obtain for $\lambda = \lambda_1$,
 \beqnn
\mbf{P}\{\tau_a^-<\infty\}
\ge
\mbf{P}\{\tau_a^-\le\tau_b^+\wedge t_0\}
\ge [g(X_0)-g(b)-\e^{-\lambda^{3/2}d_0}]/g(a)>0.
 \eeqnn

To show that $\mathbf{P} \{\tau_b^+ < \infty\} > 0$, let $m$ be a continuous and non-decreasing function with $m'(x)=0$ for $x>b+2$ and $m(x)=x$ for $0<x<b+1$.
Define
$h(x)=\e^{\lambda m(x)}$ for $\lambda,x>0$.
Then $h$ is a bounded function, and
for all $a\le x\le b$,
 \beqnn
\int_0^\infty D_zh(x)\mu(\dd z)
 \ar\ge\ar
\int_0^1 D_zh(x)\mu(\dd z)
-\lambda h(x)\int_1^\infty z\mu(\dd z) \cr
 \ar=\ar
h(x)\int_0^1 [\e^{\lambda z}-1-\lambda z]\mu(\dd z)
-\lambda h(x)\int_1^\infty z\mu(\dd z).
 \eeqnn
Using Taylor expansion of $\mathrm{e}^{\lambda z}$ around $z=0$ up to the second order,
there is a constant $\lambda_2>0$ such that for all $\lambda>\lambda_2$,
 \beqnn
[\lambda h(x)]^{-1}\int_0^\infty D_zh(x)\mu(\dd z)
\ge \lambda\int_0^1z^2\mu(\dd z)
-\int_1^\infty z\mu(\dd z)\ge \lambda^{3/4}.
 \eeqnn
Choose $\lambda_3 > \lambda_2$ such that for all $\lambda > \lambda_3$, we have $\e^{\lambda X_0} - \e^{\lambda a} - \e^{\lambda b-\lambda^{3/2}b_1((0,t_0])} > 0$ and
 \beqnn
2^{-1}\lambda  \inf_{x\in[a,b]}\gamma_1(x)
\ge
\sup_{a\le x\le b,s\in(0,t_0]}|\gamma_0(s,x)|/b_1(s),\qquad
~\sqrt{\lambda}\inf_{x\in[a,b]}
\gamma_1(x)\ge2.
 \eeqnn
Then, for $a \le x \le b$, $0 < s \le t_0$, and $\lambda > \lambda_3$,
 \beqlb\label{3.13bb}
\mathcal{L}_sh(x)
 \ar\ge\ar
\lambda h(x)
[\gamma_0(s,x)+\lambda\gamma_1(x)b_1(s)
+\lambda^{1/2}\gamma_2(x)b_2(s)] \cr
\ar\ge \ar
\lambda h(s)b_1(s)[-|\gamma_0(s,x)|/b_1(s)+2^{-1}\lambda \gamma_1(x)]+2^{-1}\lambda^2 \gamma_1(x)b_1(s)h(x) \cr
 \ar \ge \ar
\lambda^{3/2}b_1(s)h(x).
 \eeqlb
Applying a similar argument again, we find
 \beqnn
h(X_0)
 \ar\le\ar
\mbf{E}\Big[h(X_{\gamma_{a,b}})
\e^{-\lambda^{3/2}b_1((0,\gamma_{a,b}])}
\Big(1_{\{\tau_a^-\le\tau_b^+\wedge t_0\}}
+1_{\{\tau_b^+\le\tau_a^-\wedge t_0\}}
+1_{\{t_0<\tau_a^-\wedge \tau_b^+\}}\Big)\Big] \cr
 \ar\le\ar
h(a)+\mbf{P}\{\tau_b^+\le\tau_a^-\wedge t_0\}
+
h(b)\e^{-\lambda^{3/2}b_1((0,t_0])},
 \eeqnn
which implies $\mbf{P}\{\tau_b^+<\infty\}>0$.
\qed

\section{Contraction in a weighted-total variation distance}

In this section, we always assume that the SDE \eqref{ii1.5} possesses a unique strong solution $X$, the solution is non-explosive and exhibits the Markov property.

Given a nonnegative Borel function $V$ on $\mathbb{R}_+$, we denote by ${\mathcal{P}}_V(\mathbb{R}_+)$ the set of Borel probability measures $\gamma$ on $\mathbb{R}_+$ such that
$$
\int_{\mathbb{R}_+}V(x)\gamma(\dd x)<\infty.
$$
Let $W_V$ be the {\it $V$-weighted total variation distance} on ${\mathcal{P}}_V(\mathbb{R}_+)$ defined by
$$
W_V(\gamma,\eta)=\int_{\mathbb{R}_+}[1+V(x)]|\gamma-\eta|(\dd x),\quad \gamma,\eta\in {\mathcal{P}}_V(\mathbb{R}_+),
$$
where $|\cdot|$ denotes the norm of total variation. We shall see that $W_V$ is actually the Wasserstein distance determined by the metric
\beqlb\label{metric of dv}
d_V(x,y)=[2+V(x)+V(y)]\mathbf{1}_{\{x\neq y\}},\quad x,y\in\mathbb{R}_+.
\eeqlb
More precisely, we have
\beqlb\label{wasserstein dist}
W_V(\gamma,\eta)=\inf_{\pi\in{\mathcal{C}}(\gamma,\eta)}\int_{\mathbb{R}_+^2}d_V(x,y)\pi(\dd x,\dd y),
\eeqlb
where ${\mathcal{C}}(\gamma,\eta)$ is the collection of all probability measures on $\mathbb{R}_+^2$ with marginals $\gamma$ and $\eta$;  see \cite[Lemma 5.3]{LLWZ25}.  This distance has been used widely in the theory of ergodicity; see  \cite{EGZ19, HM11, LLWZ25,LMW21}. In particular, if $V\equiv0$, then $W_V$ reduces to the {\it total variation distance}. The other two frequently used weighted functions are given by
\beqnn
V_1(x)=x\quad \text{and} \quad V_{\log}(x)=\log(1+x),\quad x\ge0.
\eeqnn

\subsection{Main result}

    Let $\delta_x$ denote the Dirac measure concentrated on $x$. Given $\sigma$-finite measures $\nu_1$ and $\nu_2$ on $\mbb{R}_+$, we write $\nu_1\land \nu_2 := \nu_1 - (\nu_1-\nu_2)^+ = \nu_2 - (\nu_2-\nu_1)^+$, where the subscript ``$+$'' stands for the upper variation of the signed measure in its Jordan decomposition. We first present conditions on  SDE (\ref{ii1.5}).

\bcondition\label{cond noise}
{\bf (A) (Drift condition)} There exists $l_0>0$, such that for any $l>l_0\vee1$, there is $k_0(l)>0$ such that the following holds for any $s\ge0$,
$$
\frac{(\gamma_0(s,x)-\gamma_0(s,y))(x-y)}{|x-y|}\le k_0(l)|x-y|b_0(s),\quad  |x-y|\le l
$$
and
$$
(\gamma_0(s,x)-\gamma_0(s,y))(x-y)\le0,\quad |x-y|>l,
$$
where $b_0$ is some positive and measurable function, and $\gamma_0(s,0)\ge b_0(s)\gamma_0$ for some $\gamma_0>0$.

{\bf (B)} One of the following two assumptions holds,

(B1) {\bf (Diffusion condition)}

(B1-i)~
$\gamma_1(\cdot)\not\equiv0$;\quad (B1-ii) $\sigma b_0(s)\le b_1(s)$ holds for some $\sigma>0$ and all $s\ge0$;

(B1-iii)~for any $l>0$, there exists $k_1(l)>0$ such that
$$
\gamma_1(x)+\gamma_1(y)\ge k_1(l)|x-y|,\quad |x-y|\le l.
$$

(B2) {\bf (Jump condition)}

(B2-i)~$\gamma_1(\cdot)\equiv0$;\quad (B2-ii)~$\int_0^1z\,\mu(\dd z)=\infty$, and there are constants $c_0, \kappa_0>0$ such that
$$
(\mu\wedge(\delta_x\ast\mu))(\mathbb{R}_+)\ge\kappa_0,\quad |x|\le c_0.
$$

(B2-iii) $\bar{\sigma} b_0(s)\ge b_2(s)\ge\underline{\sigma} b_0(s)$
holds for some $\bar{\sigma}\ge\underline{\sigma}>0$ and for all $s\ge0$;

(B2-iv)~ $\lim_{x\rightarrow0}\gamma_2(x)=\gamma_2(0)=0$, $x\mapsto\gamma_2(x)$ is increasing, and for any $l>0$, there exists $k_2(l)>0$ such that
$$
|\gamma_2(x)-\gamma_2(y)|\ge k_2(l)|x-y|,\quad |x-y|\le l.
$$

{\bf (C) (Lyapunov type condition)} There exist a nonnegative function $V(x)\rightarrow\infty$ as $x\rightarrow\infty$ and $\lambda_1,\lambda_2>0$ such that for any $x\ge0$ and $s\ge0$,
$$
\mathcal{L}_sV(x)\le [-\lambda_1V(x)+\lambda_2]b_0(s),
$$
where $\mathcal{L}_s$ is the operator defined by \eqref{generator of X}.
\econdition

For a given $s\ge0$, consider the time-inhomogeneous process $(X_t)_{t\ge s}$ and let $P_{s,t}(x,\cdot)$ denote its transition probabilities for any $t\ge s$ and $x\in\mathbb{R}_+$. The main result of this section is the following theorem.

\begin{theorem}\label{exponential ergod}
Suppose Condition \ref{cond noise} is satisfied.
Then there exist constants $C_0,C_1>0$ such that for any $t>s\ge0$ and $x,y>0$,
$$
W_V(P_{s,t}(x,\cdot),P_{s,t}(y,\cdot))\le C_0 d_V(x,y)\mathrm{exp}\left(-C_1\int_s^tb_0(r)\,\dd r\right),
$$
where $d_V$ and $W_V$ denote the metric defined in \eqref{metric of dv} and the Wasserstein distance introduced in \eqref{wasserstein dist}, respectively. In particular, if
$$
\limsup\limits_{t\rightarrow\infty}\left( t^{-1}\int_\cdot^tb_0(r)\,\dd r\right)>0,
$$
the process is exponentially contractive in the $W_V$-distance.
\end{theorem}

We make some comments on Condition \ref{cond noise} and Theorem \ref{exponential ergod}.
\begin{remark} \rm
(1) Theorem \ref{exponential ergod} provides two major extensions over the main results in \cite{LLWZ25}. First, the SDE \eqref{ii1.5} is time-inhomogeneous. When $b_1$ and $  b_2$ are nonzero constants and $\gamma_0(s,x)=\gamma_0(0,x)$ for all $s\ge0$ and $x\ge0$, it reduces to a time-homogeneous SDE, as studied in \cite{LYZh}. Second, the functions $\gamma_1$ and $\gamma_2$ in \eqref{ii1.5} are allowed to be general nonlinear functions. In particular, if both $\gamma_1$ and $\gamma_2$ are linear, the model \eqref{ii1.5} becomes the classical competitive branching process discussed in \cite{LLWZ25,P16} and other works. Theorem \ref{exponential ergod} shows that,  under comparable conditions on $b_0,b_1$ and $b_2$, the study of the ergodic rate of the process $(X_t)_{t\ge s}$ yields results that go beyond exponential contraction. Moreover, the function
$b_0$ is not required to be continuous.

(2) Condition \ref{cond noise} (A) is an assumption on the drift term, under which the time and space variables are separated. In the spatial variable, we assume a local Lipschitz condition near the origin and monotonic decay at infinity, while  the temporal variable is controlled by a measurable function $b_0$. This assumption on the spatial variable is standard in ergodic theory for the time-homogeneous case.
In the absence of Condition \ref{cond noise} (C), a stronger dissipativity condition must be imposed; we refer to \cite{APS19, APS22, LW20, LMW21,LW19} for details. The time-homogeneous version of Condition \ref{cond noise} (C) is a standard assumption in the ergodic theory of Markov processes; see, for example, \cite{DMT95, MT93}. The connection between the Lyapunov condition and dissipativity can be found in \cite{EGZ19, LLWZ25}.

(3) Condition \ref{cond noise} (B) is an assumption on the stochastic integral terms. We do not require the presence of both the diffusion and jump terms simultaneously. The requirements on the spatial variable can be directly adapted from \cite[Theorem 3.1]{LW20}. In particular, since the jump term corresponds to a nonlocal operator and involves more elaborate estimates, stronger assumptions are imposed  compared to the diffusion term.
\end{remark}

To illustrate the contribution of Theorem \ref{exponential ergod}, we first provide a brief overview of contractive and ergodic properties of Markov processes, including both time-homogeneous and time-inhomogeneous cases.

In the time-homogeneous setting, the model's stability is often characterized by the contraction property of the semigroup, meaning that
$W_{V}(P_t(x,\cdot),P_t(y,\cdot))$ decays to zero as $t\rightarrow\infty$ for all $x,y\ge0$. If this decay is exponential, standard arguments (see, e.g., \cite[p.602]{FJKR23} or \cite[p.381]{LLWZ25}) imply the existence of a unique stationary distribution $\gamma\in{\mathcal{P}}_V(\mathbb{R}_+)$ such that $W_V(P_t(x,\cdot),\gamma)$ also
decays exponentially for all $x\ge0$. In this case, we say the process (or the semigroup $(P_t)_{t\ge0}$) is {\it exponentially ergodic under the $W_V$-distance}.

In contrast, results for time-inhomogeneous models are relatively scarce in the literature. Furthermore, in this setting, stationary  distributions (if they exist) typically form a family indexed by time, rather than a single distribution. This fact renders standard approaches for homogeneous processes (as in \cite[p.602]{FJKR23} or \cite[p.381]{LLWZ25}) inapplicable. In this work, we develop a framework for studying semigroup contraction in the time-inhomogeneous case, which not only incorporates existing results for the time-homogeneous setting but also reveals several new phenomena specific to the time-inhomogeneous case.

Two mainstream approaches for proving ergodicity are the Foster-Lyapunov criterion (see \cite{MT93} for example) and the coupling method (see \cite{chen04} for instance). The former typically requires verifying the irreducibility of a skeleton chain, while the latter involves constructing a suitable coupling process, defining an associated coupling operator (see the next subsection), and selecting an appropriate test function. A key advantage of the coupling method is that it often yields explicit convergence rates. Our  proof of Theorem \ref{exponential ergod} is based on the coupling approach. To the best of our knowledge, applications of coupling methods to time-inhomogeneous Markov processes remain relatively limited, and our work provides a novel contribution in this direction.

To illustrate the applicability of our results, we present two examples based on Theorems \ref{strong solution of X},
\ref{t1.3} and \ref{exponential ergod}.

\bexample\label{e4.1}
Suppose that $\gamma_1(x)=x$ and $\gamma_2$ satisfies
Condition \ref{cond for solution} (ii), $b_1$ and $b_2$ are nonnegative, measurable functions and bounded on any closed interval,
and that
there are constants $a_0,\theta\ge0$ and $a_1>0$ such that $\gamma_0(s,x)=a_1(1-a_0x)s^\theta$ for all $x,s\ge0$.
Then there is a pathwise unique solution to \eqref{ii1.5},
and the solution is
non-explosive and has the Markov property.
Moreover, the solution is exponentially contractive in the $W_V$-distance if
$\inf_{s>0}b_1(s)s^{-\theta}>0$.
\eexample

The mean field theory simplifies complex systems by having each component interact with a shared average field. This approach is now widely applied in statistical physics, financial mathematics, and social sciences. The following example indicates that \eqref{ii1.5} can be used to study
the long-time behavior of solution to a mean field SDE.

\bexample\label{e4.2}
We consider the following killed
mean field stochastic equation:
 \beqlb\label{10.1}
Z_t=Z_0+a\int_0^t\mbf{E}[Z_{s\wedge\tilde{\tau}_0}](1-bZ_s)\dd s
+\int_0^t\int_0^{Z_s} \tilde{K}(\dd s,\dd u),\quad t< \tilde{\tau}_0,
 \eeqlb
where $\tilde{\tau}_0$ is defined analogously to  $\tau_0$ by replacing  $(X_t)_{t\ge0}$
with $(Z_t)_{t\ge0}$ in the definition given at the end of Section 1, and
$\tilde{K}(\dd s, \dd u)$ is a Gaussian white noise with
density $2\tilde{a}\tilde{h}(s)\dd u$ and $\tilde{a}\le a$, where $\tilde{h}(s):=\mbf{E}[Z_{s\wedge\tilde{\tau}_0}]$.
Then the process $(Z_t)_{t\ge0}$ is a solution to \eqref{10.1} with $\mbf{E}[Z_{t\wedge\tilde{\tau}_0}]<\infty$ for all $t\ge0$
if and only if it is a solution to \eqref{ii1.5}
with coefficients given by
 \beqlb\label{10.1b}
\gamma_0(s,x)=a(1-bx)h(s),~\gamma_1(x)=x,~
\gamma_2(x)=0,~b_1(s)=\tilde{a}h(s),\quad x,s\ge0,
 \eeqlb
and
 \beqlb\label{10.1c}
h(s)=[(Z_0^{-1}-b)\e^{-at}+b]^{-1},\qquad
s\ge0.
 \eeqlb
Moreover,  the SDE \eqref{10.1} admits a pathwise unique solution $(Z_t)_{t\ge0}$
that is nonnegative, non-extinctive, non-explosive, Markov, and contractive in the $W_{V}$-distance. This solution satisfies the mean field stochastic equation:
 \beqlb\label{10.1a}
Z_t=Z_0+a\int_0^t\mbf{E}[Z_s](1-bZ_s)\dd s+\int_0^t\int_0^{Z_s} K(\dd s,\dd u),
 \eeqlb
where $Z_0,a,b>0$ and $K(\dd s, \dd u)$ is a Gaussian white noise with
density $2\tilde{a}\mbf{E}[Z_s]\dd s\dd u$.
\eexample
\proof
We first assume that $(Z_t)_{t\ge0}$ satisfies \eqref{10.1} and $\mbf{E}[Z_{t\wedge\tilde{\tau}_0}]<\infty$ for all $t\ge0$.  The stopping times $\tilde{\tau}_{\infty}$, $\tilde{\tau}_{n^{-1}}^-$, and $\tilde{\tau}_n^+$ are defined analogously to $\tau_\infty$, $\tau_{n^{-1}}^-$, and $\tau_{n}^+$, respectively, by replacing the process $(X_t)_{t\ge0}$ with $(Z_t)_{t\ge0}$ in the definitions provided at the end of Section 1.
Let $\tilde{\tau}_n:=\tilde{\tau}_{n^{-1}}^-\wedge \tilde{\tau}_{n}^+$.
Then $\tilde{\tau}_\infty=\infty$ almost surely.
From \eqref{10.1},  we have
 \beqlb\label{10.2}
\mbf{E}[Z_{t\wedge\tilde{\tau}_n}]
=Z_0+a
\mbf{E}\Big[\int_0^{t\wedge\tilde{\tau}_n}
\mbf{E}[Z_{s\wedge\tilde{\tau}_0}](1-bZ_s)\dd s\Big]
\le
Z_0+a\int_0^t\mbf{E}[Z_{s\wedge\tilde{\tau}_0}]\dd s.
 \eeqlb
Applying Fatou's lemma yields
 \beqnn
\mbf{E}[Z_{t\wedge\tilde{\tau}_0}]
\le
\liminf_{n\to\infty}\mbf{E}[Z_{t\wedge\tilde{\tau}_n}]
\le
\mbf{E}[Z_{t\wedge\tilde{\tau}_0}]\le
Z_0+a\int_0^t\mbf{E}[Z_{s\wedge\tilde{\tau}_0}]\dd s.
 \eeqnn
By Gronwall's inequality, $\tilde{h}(t)=\mbf{E}[Z_{t\wedge\tilde{\tau}_0}]\le
Z_0\e^{at}$. Moreover,
$\mbf{E}[Z_{t\wedge\tilde{\tau}_n}]\le Z_0+Z_0\int_0^t\e^{as}\dd s$.
Hence, $(Z_t)_{t\ge0}$ is a solution to \eqref{ii1.5} with  coefficients given by \eqref{10.1b} and the function $h$ replaced by $\tilde{h}$.
By Theorem \ref{t1.3},  the process $(Z_t)_{t\ge0}$ is non-extinctive,  i.e.,
$\tilde{\tau}_0=\infty$ almost surely.
Furthermore, $\lim_{n\to\infty}\tilde{\tau}_n=\infty$
almost surely.  Applying It\^o's formula to \eqref{10.1},  we obtain
 \beqnn
Z_t^2\le Z_0^2+2(a+\tilde{a})\int_0^tZ_s
\mbf{E}[Z_s]\dd s
+\mbox{martingale}, \qquad t\le \tau_n,
 \eeqnn
which implies
 \beqnn
\mbf{E}[Z_{t\wedge\tau_n}^2]\le Z_0^2+2(a+\tilde{a})\int_0^t\mbf{E}[Z_{s\wedge\tau_n}]
\mbf{E}[Z_s]\dd s
\le C(t),
 \eeqnn
where $C(t)>0$ is independent of $n$.
Using \eqref{10.2}, we conclude
 \beqnn
\mbf{E}[Z_t]
 \ar=\ar
\lim_{n\to\infty}
\mbf{E}[Z_{t\wedge\tilde{\tau}_n}]
=Z_0+a
\lim_{n\to\infty}\mbf{E}\Big[\int_0^{t\wedge\tilde{\tau}_n}
\mbf{E}[Z_s](1-bZ_s)\dd s\Big] \cr
 \ar=\ar
Z_0+a
\mbf{E}\Big[\int_0^t
\mbf{E}[Z_s](1-bZ_s)\dd s\Big]
=Z_0+a
\int_0^t
\mbf{E}[Z_s](1-b\mbf{E}[Z_s])\dd s,
 \eeqnn
which shows that
$\mbf{E}[Z_s]=h(s)$ given by \eqref{10.1c}.
Therefore,  $(Z_t)_{t\ge0}$ satisfies \eqref{10.1a} and
is a solution to \eqref{ii1.5} with coefficients given by \eqref{10.1b} and \eqref{10.1c}.

Conversely, suppose that  $(Z_t)_{t\ge0}$
is a solution to \eqref{ii1.5} with coefficients given by   \eqref{10.1b} and \eqref{10.1c}.
By Theorem \ref{t1.3},  $(Z_t)_{t\ge0}$ is non-extinctive and
non-explosive.
Moreover, $\mbf{E}[Z_s]=h(s)$, so
$(Z_t)_{t\ge0}$ satisfies both  \eqref{10.1} and \eqref{10.1a}.
The remaining assertions follows from
Theorems \ref{strong solution of X},
\ref{t1.3} and \ref{exponential ergod}.
\qed

\subsection{The coupling operator and the coupling process}

Let $X$ denote the strong solution to \eqref{ii1.5}. A strong Markov process $(X,Y)$ on $\mathbb{R}_+^2$ is called a {\it Markov coupling} of $X$ if the marginal process $Y$ shares the same transition semigroup as $X$. The Markov coupling process $(X,Y)$ satisfies $X_{T_s+t}= Y_{T_s+t}$ for all $t \ge 0$, where  $$T_s=\inf\{t\ge s: X_t = Y_t\}$$  is called the  {\it coupling time} for any $s\ge0$.

To construct a Markov coupling for the process $X$ defined by \eqref{ii1.5}, we begin
by construction of a coupling operator for its generator $({\mathcal{L}}_s)_{s\ge0}$
given in \eqref{generator of X}.  Let $({\tilde{\mathcal{L}}}_s)_{s\ge0}$ denote the infinitesimal generator of the Markov coupling process $(X,Y).$ This operator satisfies the marginal property: for any $f_1, f_2 \in C^2(\mbb{R}_+)$ and $s\ge0$,
\beqlb\label{connect_coupgen_and_gen}
{\tilde{\mathcal{L}}}_sh(x, y) = {\mathcal{L}}_sf_1(x) + {\mathcal{L}}_sf_2(y),
\eeqlb
where $h(x, y) = f_1(x) + f_2(y)$ for any $x, y \in \mbb{R}_+,$ and $({\mathcal{L}}_s)_{s\ge0}$ is given by \eqref{generator of X}. We call $({\tilde{\mathcal{L}}}_s)_{s\ge0}$ is a {\it coupling operator} of $({\mathcal{L}}_s)_{s\ge0}$.

We combine the coupling by reflection for space-time white noises with
basic coupling for Poisson random measures. The coupling by reflection for space-time white noises means using $\{-W(\dd s,\dd u): s,u\ge0\}$ (viewed as a reflection of $\{W(\dd s,\dd u): s,u\ge0\}$) before the two marginal processes meet. To explain the basic coupling for Poisson random measures from the coupling operator perspective, we begin with a coupling for the L\'evy measure $\mu$. We introduce the notation
$$
x\rightarrow x+z,\quad \mu(\dd z)
$$
for a transition from a point $x$ to the point $x+z$ with jump intensity $\mu(\dd z)$. The essential idea
of basic coupling is to make the two marginal processes jump to the same point with the biggest possible rate, where the biggest jump rate is the maximal common part of the jump intensities. In the L\'evy setting, this takes the form
$$
\mu_{y-x}(\dd z):=[\mu\wedge(\delta_{y-x}\ast\mu)](\dd z),
$$
where $x$ and $y$ represent the positions of the two marginal processes before the jump. Then basic coupling of $N(\dd s, \dd z, \dd u)$  is then constructed through the following transitions: for any $x,y\in\mathbb{R}_+$,
\beqnn
(x,y)\rightarrow
\begin{cases}
	(x+z,y+z+x-y), & \frac{1}{2}[\gamma_2(x)\wedge\gamma_2(y)]\mu_{y-x}(\dd z),\\
    (x+z,y+z+y-x), & \frac{1}{2}[\gamma_2(x)\wedge\gamma_2(y)]\mu_{x-y}(\dd z),\\
    (x+z,y+z), & [\gamma_2(x)\wedge\gamma_2(y)]
    \left(\mu-\frac{1}{2}\mu_{y-x}-\frac{1}{2}\mu_{x-y}\right)(\dd z),\\
    (x+z,y), & [\gamma_2(x)-\gamma_2(y)]^+\mu(\dd z),\\
    (x,y+z), & [\gamma_2(x)-\gamma_2(y)]^-\mu(\dd z).
\end{cases}
\eeqnn
Based on this framework, the operator corresponding to the reflection coupling and the basic coupling is defined as follows: for any $F\in C^2(\mathbb{R}_+^2)$, $s\ge0$,
\begin{equation}\label{coupling generator}
\tilde{\mathcal{L}}_sF(x,y)=\tilde{\mathcal{L}}_s^{d}F(x,y)+b_1(s)\tilde{\mathcal{L}}^{c}F(x,y)+b_2(s)\tilde{\mathcal{L}}^{j}F(x,y),
\end{equation}
where
\beqlb\label{coupling generator_1}
\tilde{\mathcal{L}}_s^dF(x,y)=\gamma_0(s,x)F'_x(x,y)
+\gamma_0(s,y)F'_y(x,y),
\eeqlb
\beqlb\label{coupling generator_2}
\tilde{\mathcal{L}}^{c}F(x,y)=\gamma_1(x)F''_{xx}(x,y)+\gamma_1(y)F''_{yy}(x,y)
-2\gamma_1(y)F''_{xy}(x,y)
\eeqlb
and
\beqlb\label{coupling generator_3}
\tilde{\mathcal{L}}^{j}F(x,y)
\ar=\ar
[\gamma_2(x)\wedge\gamma_2(y)]\int_0^\infty D_{(z,z+x-y)}F(x,y)
\,(\mu_{y-x}/2)(\dd z)\cr
\ar\ar\,+
[\gamma_2(x)\wedge\gamma_2(y)]\int_0^\infty D_{(z,z+y-x)}F(x,y)
\,(\mu_{x-y}/2)(\dd z)\cr
\ar\ar\,+
[\gamma_2(x)\wedge\gamma_2(y)]\int_0^\infty D_{(z,z)}F(x,y)
\left(\mu-\mu_{y-x}/2-\mu_{x-y}/2\right)(\dd z)\cr
\ar\ar\,+[\gamma_2(x)-\gamma_2(y)]^+\int_0^\infty D_{(z,0)}F(x,y)
\,\mu(\dd z)\cr
\ar\ar\,+
[\gamma_2(x)-\gamma_2(y)]^-\int_0^\infty
D_{(0,z)}F(x,y)\,\mu(\dd z),
\eeqlb
here and in what follows, we denote $F'_x(x,y)=\frac{\partial F}{\partial x}(x,y)$ and $F''_{xx}(x,y)=\frac{\partial^2 F}{\partial x^2}(x,y)$ with similar notations used for other partial derivatives. The operator $\Delta_{(z_1,z_2)}$ is defined at the end of Section 1.2. One can verify that $(\tilde{\mathcal{L}}_s)_{s\ge0}$ defined by \eqref{coupling generator} is indeed a coupling generator with respect to $({\mathcal{L}}_s)_{s\ge0}$
defined in \eqref{generator of X}; see \cite[Subsection 2.1]{LW19} for similar discussions in the time-homogeneous coefficient case.

Consider the following system of SDEs for $t>s\ge0$,
\beqlb\label{coupling process}
X_t\ar=\ar X_s+\int_s^t\gamma_0(r-,X_{r-})\dd r+\int_s^t\int_0^{\gamma_1(X_{r-})}\,W(\dd r,\dd  u)+\int_s^t\int_0^\infty\int_0^{\gamma_2(X_{r-})} z\,\tilde{N}(\dd r,\dd z,\dd u),\cr
Y_t\ar=\ar Y_s+\int_s^t\gamma_0(r-,Y_{r-})\dd r
-\int_s^{t\wedge T_s}\int_0^{\gamma_1(Y_{r-})}\,W(\dd r,\dd u)+
\int_{t\wedge T_s}^t\int_0^{\gamma_1(Y_{r-} )}\,W(\dd r,\dd u)\cr
\ar\ar\,
+\int_s^t\int_0^\infty\int_0^{\frac{1}{2}\gamma_2(Y_{r-})\rho(-U_{r-},z)} [z+U_{r-}]\,\tilde{N}(\dd r,\dd z,\dd u)\cr
\ar\ar\,
+\int_s^t\int_0^\infty\int_{\frac{1}{2}\gamma_2(Y_{r-})\rho(-U_{r-},z)}
^{\frac{1}{2}\gamma_2(Y_{r-})(\rho(-U_{r-},z)+\rho(U_{r-},z))} [z-U_{r-}]\,\tilde{N}(\dd r,\dd z,\dd u)\cr
\ar\ar\,
+\int_s^t\int_0^\infty\int_
{\frac{1}{2}\gamma_2(Y_{r-})(\rho(-U_{r-},z)+\rho(U_{r-},z))}^{\gamma_2(Y_{r-})} z\,\tilde{N}(\dd r,\dd z,\dd u),
\eeqlb
where $T_s=\inf\{t\ge s:  X_t=Y_t\}$, $U_r=X_r-Y_r$ for each $r\ge0$, and the control function
$$
\rho(x,z)=\frac{\mu_x(\dd z)}{\mu(\dd z)}\in[0,1],\quad x\in\mathbb{R},\quad z\in\mathbb{R}_+
$$
with $\rho(0,z)=1$ by convention. From \cite[Corollary
A.2 and Remark 2.1]{LW19}, we have
\beqlb\label{mu_z,mu_-z}
\mu_x(\mathbb{R}_+)=\mu_{-x}(\mathbb{R}_+)<\infty, \quad \text{for~any}~~
x\in\mathbb{R}.
\eeqlb

\begin{proposition}
For any $(X_s, Y_s)=(x,y)\in\mathbb{R}_+^2$, the system \eqref{coupling process} is well-defined, and admits a unique strong solution whose infinitesimal generator coincides with the coupling operator
defined in \eqref{coupling generator}.
\end{proposition}

\proof  It suffices to establish the existence and uniqueness of a strong solution $Y$  given $X$.} For simplicity, we only consider $s=0$. We show that the sample paths of
$Y$ defined in \eqref{coupling process} can be constructed through iterative modifications of the strong solution to the following SDE:
\beqlb\label{process z}
Z_t=y+\int_0^t\gamma_0(s-,Z_{s-})\dd s+\int_0^t\int_0^{\gamma_1(Z_{s-})}\,W^*(\dd s,\dd  u)+\int_0^t\int_0^\infty\int_0^{\gamma_2(Z_{s-})} z\,\tilde{N}(\dd s,\dd z,\dd u),
\eeqlb
where
\beqnn
W^*(\dd s,\dd u)
=-\mathbf{1}_{\{s\le T_0\}}W(\dd s,\dd u)+\mathbf{1}_{\{s>T_0\}}W(\dd s,\dd u).
\eeqnn
One can verify that $\{W^*(\dd s,\dd u): s\ge0, u\ge0\}$ remains a Gaussian white noise on $(0,\infty)^2$ with intensity $2b_1(s-)\mathrm{d}s\mathrm{d}u$. Since the driving Poisson random measure for \eqref{ii1.5} and \eqref{process z} is identical, the existence of a strong solution $Z$ to \eqref{process z} follows from the pathwise unique
strong solution to \eqref{ii1.5}.

We now claim that the process $Y$ defined in \eqref{coupling process} coincides with the solution of the following SDE:
\beqlb\label{process Y2}
Y_t\ar=\ar y+\int_0^t\gamma_0(s-,Y_{s-})\dd s
+\int_0^t\int_0^{\gamma_1(Y_{s-})}\,W^*(\dd s,\dd u)\cr
\ar\ar\,
+\int_0^t\int_0^\infty\int_0^{\gamma_2(Y_{s-})}z\,\tilde{N}(\dd s,\dd z,\dd u)\cr
\ar\ar\,
+\int_0^tU_{s-}\int_0^\infty\int_0^{\frac{1}{2}\gamma_2(Y_{s-})\rho(-U_{s-},z)}
\,N(\dd s,\dd z,\dd u)\cr
\ar\ar\,
-\int_0^tU_{s-}\int_0^\infty
\int_
{\frac{1}{2}\gamma_2(Y_{s-})\rho(-U_{s-},z)}
^{\frac{1}{2}\gamma_2(Y_{s-})[\rho(-U_{s-},z)+\rho(U_{s-},z)]} \,N(\dd s,\dd z,\dd u).
\eeqlb

\eqref{process Y2} follows directly from \eqref{mu_z,mu_-z} and the identity that for any $x,y\in\mathbb{R}_+$ with $x\neq y$,
\beqnn
\int_0^\infty\int_0^{\frac{1}{2}\gamma_2(y)\rho(-(x-y),z)}\dd u\mu(\dd z)&=&\frac{1}{2}\gamma_2(y)\mu_{-(x-y)}(\mathbb{R}_+)=
\frac{1}{2}\gamma_2(y)\mu_{(x-y)}(\mathbb{R}_+)\\
&=&\int_0^\infty\int_{\frac{1}{2}\gamma_2(y)\rho(-(x-y),z)}
^{\frac{1}{2}\gamma_2(y)[\rho(-(x-y),z)+\rho((x-y),z)]}\dd u\mu(\dd z).
\eeqnn
Furthermore, adapting the proof technique of  \cite[Proposition 2.2]{LW20} by defining appropriate stopping times and constructing the trajectories of $Y$ incrementally using the processes $X$ and $Z$ from \eqref{ii1.5} and \eqref{process z}, we conclude that the unique strong solution $Y$ to the SDE \eqref{process Y2}  can be determined globally.

Applying It\^{o}'s formula to the system \eqref{coupling process} shows that the infinitesimal generator of the process $(X,Y)$ in \eqref{coupling process} coincides with the coupling operator defined in \eqref{coupling generator}.\qed

For our main results, the crucial property of this coupling is its ability to bring the two processes together. The following estimates on the generator, derived for specific test functions, will be essential in later subsections to formalize this property.

Suppose that
\beqlb\label{test function1}
F(x,y)=f(|x-y|)\mathbf{1}_{\{x\neq y\}},
\eeqlb
where $f$ is a nonnegative, nondecreasing, bounded, and concave function in $C^2(\mathbb{R}_+)$.
By Taylor expansion, we have
\beqnn
D_{(z,z+x-y)}F(x,y)
=\frac{1}{2}\left(z^2\frac{\partial^2}{\partial x^2}
+(z+x-y)^2\frac{\partial^2}{\partial y^2}\right)F\Big(x+\theta_0z,y+\theta_0(z+x-y)\Big)
\eeqnn
for some $\theta_0\in(0,1)$. Then the first integral on the right hand side of \eqref{coupling generator_3} is nonpositive. Similar arguments apply to the remain integrals. Consequently,
\beqlb\label{coupling generator_333}
\tilde{\mathcal{L}}^{j}F(x,y)\le0.
\eeqlb

Alternatively, suppose that
$$
F(x,y)=\phi(x\vee y)\psi(|x-y|)\mathbf{1}_{\{x\neq y\}},
$$
where $\phi$ is a nonnegative, nonincreasing and bounded function in $C^2(\mbb{R}_+)$, and $\psi$ is a nonnegative, nondecreasing, bounded and concave function in $C^2(\mbb{R}_+)$. Then for $x>y$, we have $F(x,y)=\phi(x)\psi(x-y)$, and appealing to \eqref{coupling generator_3} and \eqref{mu_z,mu_-z} yields
    \beqlb\label{est of Lj1 nonsym}
    \tilde{\mathcal{L}}^jF(x,y)\ar=\ar
    \frac{1}{2}\gamma_2(y)\int_0^\infty\Big[F(x+z,x+z)-F(x,y)-z\phi'(x)\psi(x-y)\cr
    \ar\ar\,\,\,\,\qquad\quad
    +(x-y)\phi(x)\psi'(x-y)\Big]\,\mu_{y-x}(\dd z)\cr
    \ar\ar\,\,+\frac{1}{2}\gamma_2(y)\int_0^\infty\Big[F(x+z,2y-x+z)-F(x,y)-z\phi'(x)\psi(x-y)\cr
    \ar\ar\,\,\,\,\qquad\quad
    +(y-x)\phi(x)\psi'(x-y)\Big]\,\mu_{x-y}(\dd z)\cr
\ar\ar\,\,+\gamma_2(y)\int_0^\infty\Big[F(x+z,y+z)-F(x,y)-z\phi'(x)\psi(x-y)\Big]\cr
   \ar\ar\,\,\,\,\qquad\quad \times\left(\mu-\frac{1}{2}\mu_{x-y}-\frac{1}{2}\mu_{y-x}\right)(\dd z)\cr
    \ar\ar\,\,+[\gamma_2(x)-\gamma_2(y)]\int_0^\infty
    \Big[F(x+z,y)-F(x,y)-z\phi'(x)\psi(x-y)\cr
\ar\ar\,\,\,\,\qquad\quad
    -z\phi(x)\psi'(x-y)\Big]\,\mu(\dd z)\cr
    \ar=\ar \gamma_2(y)\psi(x-y)\int_0^\infty D_z\phi(x)
    \,\mu(\dd z)+\frac{1}{2}\gamma_2(y)I_1(x,y)\cr
    \ar\ar\,\,+[\gamma_2(x)-\gamma_2(y)]\int_0^\infty I_2(x,y,z)\,\mu(\dd z),
    \eeqlb
where
\beqnn
\ar\ar I_1(x,y):=[\psi(2(x-y))-\psi(x-y)]\int_0^\infty\phi(x+z)\,\mu_{x-y}(\dd z)-\psi(x-y)
\int_0^\infty\phi(x+z)\,\mu_{y-x}(\dd z),\cr
\ar\ar I_2(x,y,z)
:=\phi(x+z)\psi(x-y+z)-\phi(x)\psi(x-y)
    -z\phi'(x)\psi(x-y)-z\phi(x)\psi'(x-y).
\eeqnn
We first estimate
\beqlb\label{est of Lj1 nonsym1}
I_1(x,y)\ar=\ar-\psi(x-y)\int_0^\infty[\phi(x+z)-\phi(x)]\mu_{y-x}(\dd z)
-\psi(x-y)\phi(x)\mu_{y-x}(\mathbb{R}_+)\cr
\ar\ar\,+[\psi(2(x-y))-\psi(x-y)]\int_0^\infty\phi(x+z)\,\mu_{x-y}(\dd z)\cr
\ar\le\ar-\psi(x-y)\int_0^\infty[\phi(x+z)-\phi(x)]\mu_{y-x}(\dd z)
-\psi(x-y)\phi(x)\mu_{y-x}(\mathbb{R}_+)\cr
\ar\ar\,
+[\psi(2(x-y))-\psi(x-y)]\int_0^\infty\phi(x)\,\mu_{x-y}(\dd z)\cr
\ar=\ar-\psi(x-y)\int_0^\infty[\phi(x+z)-\phi(x)]\mu_{y-x}(\dd z)\cr
\ar\ar\,
+[\psi(2(x-y))-2\psi(x-y)]\phi(x)\,\mu_{x-y}(\mathbb{R}_+),\quad\quad x > y,
    \eeqlb
where the last equality uses \eqref{mu_z,mu_-z}. For $I_2(x,y,z)$, we have
\beqlb\label{est of Lj1 nonsym2}
     I_2(x,y,z)\ar=\ar\phi(x)
[\psi(x-y+z)-\psi(x-y)-z\psi'(x-y)]\cr
\ar\ar
+(\phi(x+z)-\phi(x))\psi(x-y+z)-z\phi'(x)\psi(x-y)\cr
\ar\le\ar\phi(x) [\psi(x-y+z)-\psi(x-y)-z\psi'(x-y)]\cr
\ar\ar
+\psi(x-y)
[\phi(x+z)-\phi(x)-z\phi'(x)],\quad\quad x > y,
\eeqlb
where the inequality follows from $\phi(x+z)-\phi(x)\le0$ and the monotonicity of $\psi$.
Combining \eqref{est of Lj1 nonsym}, \eqref{est of Lj1 nonsym1} and \eqref{est of Lj1 nonsym2} yields
    \beqlb\label{coupling generator_3333}
    \tilde{\mathcal{L}}^jF(x,y)    \ar\le\ar\gamma_2(x)\psi(x-y)\int_0^\infty
    D_z\phi(x)
    \,\mu(\dd z)\cr
    \ar\ar\,+(\gamma_2(x)-\gamma_2(y))\phi(x)
    \int_0^\infty
    D_z\psi(x-y)
   \,\mu(\dd z)\cr
    \ar\ar\,+\frac{1}{2}\gamma_2(y)\left[\psi(2(x-y))-2\psi(x-y)\right]\phi(x)\mu_{x-y}(\mathbb{R}_+)\cr
    \ar\ar\,
    -\frac{1}{2}\gamma_2(y)\psi(x-y)\int_0^\infty \Delta_z\phi(x)\,\mu_{y-x}(\dd z),\quad\quad x > y.
    \eeqlb

\subsection{Key estimates}

To obtain refined estimates for the coupling operator, we prepare by considering two cases.

\bigskip

{\bf Case(i): Condition \ref{cond noise} (A) and (B1) hold}

\bigskip

Under Condition \ref{cond noise} (A) and (B1), for any $l > l_0 \vee 1$, define
$$
\Psi_{l}(z)=-k_0(l)z+\sigma k_1(l)z^2.
$$
Choose a constant $c(l)>0$ such that for all $z\ge c(l)$,
$$
\Psi_{l}(z)\ge1,\quad \int_{c(l)}^\infty \Psi^{-1}_{l}(z)\,\dd z\le1.
$$
Define an auxiliary function
$$
f_{l}(x)=2+\int_{c(l)}^\infty\frac{1-\mathrm{e}^{-z x}}{\Psi_{l}(z)}
\,\mathrm{d}z
$$
and let
\beqlb\label{test function1}
F_{l}(x,y)=f_{l}(|x-y|)\mathbf{1}_{\{x\neq y\}}.
\eeqlb
By \eqref{coupling generator_2} and the fact that $\frac{\partial F^2_l}{\partial x^2}(x,y)=\frac{\partial F^2_l}{\partial y^2}(x,y)=-\frac{\partial F^2_l}{\partial xy}(x,y)$, we have
$$
\tilde{\mathcal{L}}^{c}F_l(x,y)=[\gamma_1(x)+3\gamma_1(y)]\frac{\partial F^2_l}{\partial x^2}(x,y).
$$
Since $f'_{l}\ge0$ and $f''_{l}\le0$, it follows from \eqref{coupling generator_3}, \eqref{test function1},  \eqref{coupling generator_333} and Condition \ref{cond noise}   (B1-ii) that for any $x\neq y$ and $s\ge0$, \eqref{coupling generator} yields
\beqlb\label{est of coup when diff}
\tilde{\mathcal{L}}_sF_{l}(x,y)\ar\le\ar [\gamma_0(s,x)-\gamma_0(s,y)]f_l'(|x-y|)+b_1(s)[\gamma_1(x)+3\gamma_1(y)]f_l''(|x-y|)\cr
\ar\le\ar [\gamma_0(s,x)-\gamma_0(s,y)]f_l'(|x-y|)+\sigma b_0(s)[\gamma_1(x)+\gamma_1(y)]f_l''(|x-y|).
\eeqlb

\begin{proposition}\label{prop10}
Suppose that Condition \ref{cond noise} (A) and (B1) hold. Then for any $l > l_0 \vee 1$, there exists a constant $\lambda(l)>0$ such that for any $s\ge0$ and $x\neq y$,
$$
\tilde{\mathcal{L}}_sF_{l}(x,y)\le-\lambda(l)b_0(s)F_{l}(x,y)\mathbf{1}_{\{0<|x-y|\le l\}}.
$$
\end{proposition}

\proof Assume $x>y$ without loss of generality (the case $x<y$ is analogous). Fix some $s\ge0$ and $l > l_0 \vee 1$. By \eqref{est of coup when diff} and Condition \ref{cond noise} (A), for $x-y>l$, we have
\beqnn
&&\tilde{\mathcal{L}}_sF_{l}(x,y)\le0.
\eeqnn
Now consider $x-y\le l$.  By Condition \ref{cond noise} (A), (B1-iii) and \eqref{est of coup when diff}, for all $x>y\ge0$ and $s\ge0$,
\beqnn
\tilde{\mathcal{L}}_sF_{l}(x,y)
\ar\le\ar k_0(l) b_0(s)(x-y)f'_l(x-y)+\sigma k_1(l) b_0(s)(x-y)f''_l(x-y)\cr
\ar=\ar k_0(l_1) b_0(s)(x-y)\int_{c(l)}^\infty z\mathrm{e}^{-z(x-y)}
\Psi_{l}(z)^{-1}\,\dd z\cr
\ar\ar\,
-\sigma k_1(l)b_0(s)(x-y)
\int_{c(l)}^\infty z^2\mathrm{e}^{-z(x-y)}
\Psi_{l}(z)^{-1}\,\dd z\cr
\ar=\ar -b_0(s)\int_{c(l)}^\infty(x-y)\mathrm{e}^{-z(x-y)}\,\mathrm{d}z=-b_0(s)\mathrm{e}^{-c(l)(x-y)}.
\eeqnn
Therefore,
$$
\tilde{\mathcal{L}}_sF_{l}(x,y)\le-b_0(s)\mathrm{e}^{-c(l)l}\mathbf{1}_{\{x-y\le l\}}.
$$
The assertion follows since
$2\le F_{l}(x,y)\le3$ whenever $x\neq y$.
\qed

\bigskip

 {\bf Case (ii): Condition \ref{cond noise} (A) and  (B2) hold}

\bigskip

Suppose Condition \ref{cond noise} (A) and  (B2) hold. For parameters $x_0\in(0,1\wedge c_0)$ (where $c_0$ is the constant from Condition \ref{cond noise} (B2-ii)) and $\theta > 1$ to be specified later, define
\beqnn
\phi(x) =
\begin{cases}
	\theta+(1-x/x_0)^3,& x \in [0,x_0),\\
	\theta, & x \in [x_0,\infty).
\end{cases}
\eeqnn

For any $l > l_0 \vee 1,$ define
    $$
    \Phi_{l}(u)=-uk_0(l)+k_2(l)\underline{\sigma}\int_0^\infty(\mathrm{e}^{-zu}-1+zu)\,\mu(\dd z),\quad u\ge0.
    $$
Now consider
    \beqnn
u^{-1}\int_0^\infty(\mathrm{e}^{-zu}-1+zu)\,\mu(\dd z)\ar=\ar \int_0^\infty z\,\mu(\dd z)\int_0^1(1-\mathrm{e}^{-zuv})\,\dd v\cr
    \ar\ge\ar \int_{u^{-1}}^1z\,\mu(\dd z)\int_{1/2}^1(1-\mathrm{e}^{-1/2})\,\dd v
    =\frac{1}{2}(1-\mathrm{e}^{-1/2})\int_{u^{-1}}^1z\,\mu(\dd z),
    \eeqnn
where the first equality follows from Taylor expansion. Hence,
    \beqnn
\Phi_{l}(u) \ge -u\Big(k_0(l)-\frac{ k_2(l)\underline{\sigma}}{2}(1-\mathrm{e}^{-1/2})\int_{u^{-1}}^1z\,\mu(\dd z)\Big).
    \eeqnn
By Condition \ref{cond noise} (B2-ii), $\int_0^1z\,\mu(\dd z)=\infty$, so there exists $c(l)>0$ such that $\Phi_{l}(c(l))>1$. Define further
    $$
    \psi_{l}(x)=2-\mathrm{e}^{-c(l)x},\quad  x\ge0
    $$
and construct the auxiliary function
\beqlb\label{test function2}
    F_{l}(x,y)=\phi(x\vee y)\psi_{l}(|x-y|)\mathbf{1}_{\{x\neq y\}}.
\eeqlb
In particular, when $x>y$,  \eqref{coupling generator_1} gives
    \beqlb\label{est of drift}
\tilde{\mathcal{L}}_s^dF_{l}(x,y)=\gamma_0(s,x)\phi'(x)\psi_{l}(x-y)
    +[\gamma_0(s,x)-\gamma_0(s,y)]\phi(x)\psi'_{l}(x-y).
    \eeqlb

\begin{proposition}\label{prop11}
Assume that Condition \ref{cond noise} (A) and (B2) hold. Then for any $l> l_0 \vee 1$, there exists a constant $\lambda(l)>0$ such that for any $s\ge0$ and $x\neq y$,
\beqlb\label{ineq est of coupling funct}
\tilde{\mathcal{L}}_sF_{l}(x,y)\le-\lambda(l) b_0(s) F_{l}(x,y)\mathbf{1}_{\{0<|x-y|\le l\}}.
\eeqlb
\end{proposition}

\proof We present the proof for the case $x > y$, the case $x < y$ can be discussed similarly.  Fix $s\ge0$ and $l > l_0 \vee 1$. By \eqref{est of drift} and Condition \ref{cond noise} (A), for $x>x_0$ and $x-y>l$, we have
$\tilde{\mathcal{L}}_s^dF_{l}(x,y)\le0$. From \eqref{coupling generator_3333}, it follows that
 \beqnn
\tilde{\mathcal{L}}^jF_l(x,y)\ar\le\ar
(\gamma_2(x)-\gamma_2(y))\phi(x)
\int_0^\infty D_z\psi_l(x-y)
\,\mu(\dd z)\cr
\ar\ar\,
+\frac{1}{2}\gamma_2(y)\left[\psi_l(2(x-y))-2\psi_l(x-y)\right]\phi(x)\mu_{x-y}(\mathbb{R}_+).
\eeqnn
Note that $D_z\psi_l(x-y)=\frac{1}{2}z^2\psi''_l(\xi)\le0$ for some $\xi\in(x-y,x-y+z)$ and
\beqlb\label{eq:4.26}
\psi_l(2(x-y))-2\psi_l(x-y)\le-1.
\eeqlb
Hence, $\tilde{\mathcal{L}}^jF_l(x,y)\le0$. Therefore,  \eqref{ineq est of coupling funct} holds for  $x>x_0$ and $|x-y|>l$.

Now consider $x>x_0$ and $x-y\le l$. From
 Condition \ref{cond noise} (A) and \eqref{est of drift}, we obtain
$$
\tilde{\mathcal{L}}_s^dF_{l}(x,y)\le c(l)k_0(l)\mathrm{e}^{-c(l)(x-y)}(x-y)\phi(x)b_0(s).
$$
Moreover, by \eqref{coupling generator_3333}, \eqref{eq:4.26} and Condition \ref{cond noise} (B2-iv), for $x>x_0$ and $x-y\le l$,
\beqnn
\tilde{\mathcal{L}}^jF_{l}(x,y)
\ar
\le\ar(\gamma_2(x)-\gamma_2(y))\phi(x)
\int_0^\infty D_z\psi_l(x-y)\,\mu(\dd z)\cr
\ar\ar\,+\frac{1}{2}\gamma_2(y)\phi(x)\mu_{x-y}(\mathbb{R}_+)[\psi_{l}(2(x-y))-2\psi_{l}(x-y)]\cr
\ar\le\ar \phi(x)\Big[-k_2(l)(x-y)\mathrm{e}^{-c(l)(x-y)}\int_0^\infty
\Big(\mathrm{e}^{-c(l)z}-1+c(l)z\Big)\,\mu(\dd z)\cr
\ar\ar\quad\quad\quad
-\mu_{x-y}(\mathbb{R}_+)\gamma_2(y)/2\Big]
< 0.
\eeqnn
Then, using  Condition \ref{cond noise} (B2-iii) and \eqref{coupling generator}, along with $\Phi_{l}(c(l))>1$, for $x>x_0$ and $x-y\le l$,
\beqlb\label{x>x0}
\tilde{\mathcal{L}}_sF_{l}(x,y)\ar=\ar \tilde{\mathcal{L}}^d_sF_{l}(x,y)+b_2(s)\tilde{\mathcal{L}}^{j}F(x,y)\cr
\ar\le\ar b_0(s)\phi(x)\left[-\Phi_{l}(c(l))\mathrm{e}^{-c(l)(x-y)}(x-y)-\frac{\underline{\sigma}}{2}\mu_{x-y}(\mbb{R}_+)\gamma_2(y)
\right]\cr
\ar
\le\ar b_0(s)\phi(x)\left[-\mathrm{e}^{-c(l)(x-y)}(x-y)
-\frac{\underline{\sigma}}{2}\mu_{x-y}(\mbb{R}_+)\gamma_2(y)
\right].
\eeqlb

For $x\le x_0$, clearly $x-y\le l$.
In this case, by \eqref{coupling generator_3333}, \eqref{est of drift}, Condition
\ref{cond
noise} (A) and (B2-iii, iv), we obtain
\beqlb\label{x<x0}
\tilde{\mathcal{L}}_sF_{l}(x,y)\ar=\ar \tilde{\mathcal{L}}^d_sF_{l}(x,y)+b_2(s)\tilde{\mathcal{L}}^{j}F(x,y)\cr
\ar
\le\ar b_0(s)\phi(x)\left[-\mathrm{e}^{-c(l)(x-y)}(x-y)
-\underline{\sigma}\mu_{x-y}(\mbb{R}_+)\gamma_2(y)/2
\right]\cr
\ar\ar
+\psi_{l}(x-y)\gamma_0(s,x)\phi'(x)\cr
\ar\ar+\bar{\sigma}b_0(s)
\gamma_2(x)\psi_{l}(x-y)\int_0^\infty D_z\phi(x)\,\mu(\dd z)\cr
\ar\ar
-\frac{1}{2}\bar{\sigma}b_0(s)\gamma_2(y)\psi_{l}(x-y)
\int_0^\infty \Delta_z\phi(x)\mu_{y-x}(\dd z).
\eeqlb
With these estimates, we proceed to prove \eqref{ineq est of coupling funct} by considering the following
three cases.

{\bf (a) $x>x_0$.} If  $x_0/2<x-y\le l$, then by \eqref{x>x0},
\beqnn
\tilde{\mathcal{L}}_sF_{l}(x,y)\le-b_0(s)\phi(x)(x-y)\mathrm{e}^{-c(l)(x-y)}
\le-\frac{x_0\mathrm{e}^{-c(l)l}}{2\psi_{l}(l)}b_0(s)F_{l}(x,y).
\eeqnn
If $x-y\le x_0/2$, then $y\ge x_0/2$, so by Condition \ref{cond noise} (B2-iv), $\gamma_2(y)\ge\gamma_2(x_0/2)$.
Since $x_0\le 1\wedge c_0$, Condition \ref{cond noise}
 (B2-ii) implies $\mu_{x-y}(\mbb{R}_+)\ge\kappa_0$.
Then from \eqref{x>x0},
$$
\tilde{\mathcal{L}}_sF_{l}(x,y)\le
-\frac{1}{2}\underline{\sigma}\kappa_0\gamma_2(x_0/2)b_0(s)\phi(x)
\le-\frac{\underline{\sigma}\kappa_0\gamma_2(x_0/2)}{2\psi_l(x_0/2)}b_0(s)F_l(x,y).
$$

{\bf (b) $x\in(0,rx_0]$}, where $r\in(0,1/2)$ will be specified later. By \eqref{x<x0}, for $x\in(0,rx_0]$,
\beqlb\label{est for x<rx0}
\tilde{\mathcal{L}}_sF_{l}(x,y)\ar\le\ar-\frac{1}{2}\underline{\sigma}b_0(s)\phi(x)
\mu_{x-y}(\mbb{R}_+)\gamma_2(y)
+\psi_{l}(x-y)\gamma_0(s,x)\phi'(x)
+\bar{\sigma}b_0(s)\psi_{l}(x-y)\cr
\ar\ar\,\,
\times\Big[
\gamma_2(x)\int_0^\infty D_z\phi(x)\,\mu(\mathrm{d}z)
-\frac{1}{2}\gamma_2(y)\int_0^\infty \Delta_z\phi(x)\mu_{y-x}(\dd z)\Big].
\eeqlb
From the definition of $\phi$, we have
\beqlb\label{est of phi}
\ar\ar 0\le -\Delta_z\phi(x)\le\mathbf{1}_{\{x\le x_0\}},\quad 0\le-\phi'(x)z\le 3zx_0^{-1}\mathbf{1}_{\{x\le x_0\}},\cr
\ar\ar
D_z\phi(x)=\phi''(\zeta)\mathbf{1}_{\{\zeta<x_0\}}z^2/2
\le 3z^2x_0^{-2}\mathbf{1}_{\{x\le x_0\}}
\eeqlb
for some $\zeta\in(x,x+z)$. Then, using   \eqref{mu_z,mu_-z}, \eqref{est of phi} and  $\phi(x)\ge\theta$, from \eqref{est for x<rx0} we obtain
\beqnn
\tilde{\mathcal{L}}_sF_{l}(x,y)\ar\le\ar-\frac{1}{2}\underline{\sigma}b_0(s)\phi(x)
\mu_{x-y}(\mbb{R}_+)\gamma_2(y)
+\psi_{l}(x-y)\gamma_0(s,x)\phi'(x)\cr
\ar\ar\,
+\bar{\sigma}b_0(s)\psi_{l}(x-y)\left[H_0\gamma_2(x)
+\frac{1}{2}\gamma_2(y)\mu_{y-x}(\mathbb{R}_+)\right]\cr
\ar\le\ar \frac{1}{2}(-\underline{\sigma}\theta+\bar{\sigma})\mu_{x-y}(\mbb{R}_+)\gamma_2(y)b_0(s)
+\psi_{l}(x-y)\gamma_0(s,x)\phi'(x)\cr
\ar\ar\,
+\bar{\sigma}H_0b_0(s)\psi_{l}(x-y)\gamma_2(x),
\eeqnn
where
$$
H_0=3x_0^{-2}\left(\int_0^1z^2\,\mu(\dd z)+
\int_1^\infty z\,\mu(\mathrm{d}z)\right).
$$
Choosing $\theta\ge \bar{\sigma}/\underline{\sigma}$ sufficiently large yields
\beqlb\label{eq:4.280}
\tilde{\mathcal{L}}_sF_{l}(x,y)\le \psi_{l}(x-y)\gamma_0(s,x)\phi'(x)+\bar{\sigma}H_0b_0(s)\psi_{l}(x-y)\gamma_2(x).
\eeqlb
By Condition \ref{cond noise} (A),
\beqlb\label{eq:4.281}
\gamma_0(s,x)\phi'(x)\ar\le\ar|\gamma_0(s,x)-\gamma_0(s,0)||\phi'(x)|+\gamma_0(s,0)\phi'(x)\cr
\ar\le\ar 3x_0^{-1}k_0(l)xb_0(s)-3\gamma_0x_0^{-1}b_0(s)(1-x/x_0)^2\cr
\ar\le\ar 3k_0(l)rb_0(s)-3\gamma_0b_0(s)/(4x_0),
\eeqlb
where the last inequality uses $x \in (0, rx_0]$ with $r < 1/2$.
By Condition \ref{cond noise} (B2-iv), \eqref{eq:4.280} and \eqref{eq:4.281}, we can choose $r\in(0,1/2)$ sufficiently small so that
\beqnn
\tilde{\mathcal{L}}_sF_{l}(x,y)\ar\le\ar b_0(s)\psi_{l}(x-y)\left[-3\gamma_0/(4x_0)+3k_0(l)r
+\bar{\sigma}H_0\gamma_2(rx_0)\right]\cr
\ar\le\ar -3\gamma_0b_0(s)\phi(x)\psi_{l}(x-y)/(8(\theta+1)x_0)\cr
\ar=\ar -3\gamma_0b_0(s)F_{l}(x-y)/(8(\theta+1)x_0).
\eeqnn

{\bf (c)  $x\in(rx_0,x_0)$.} Using an argument analogous to case (b),
from \eqref{x<x0} and \eqref{est of phi}, we obtain
\beqlb\label{eq:4.290}
\tilde{\mathcal{L}}_sF_{l}(x,y)\ar\le\ar
b_0(s)\left[-\phi(x)(x-y)\mathrm{e}^{-c(l)(x-y)}+\frac{1}{2}(-\underline{\sigma}\theta+\bar{\sigma})\mu_{x-y}(\mbb{R}_+)\gamma_2(y)\right]\cr
\ar\ar\,
+\psi_l(x-y)\gamma_0(s,x)\phi'(x)+\bar{\sigma}H_0b_0(s)\psi_l(x-y)\gamma_2(x).
\eeqlb
From \eqref{eq:4.281}, for $x\in(rx_0,x_0)$,
$$
\gamma_0(s,x)\phi'(x)\le3k_0(l)b_0(s).
$$
Since $1 \le \psi_l\le 2$ , using Condition \ref{cond noise} (B2-iv) and \eqref{eq:4.290}, we obtain
\beqlb\label{est of coup when jump (3)}
\tilde{\mathcal{L}}_sF_{l}(x,y)\le b_0(s)\left[-\phi(x)(x-y)\mathrm{e}^{-c(l)(x-y)}
+\frac{1}{2}(-\underline{\sigma}\theta+\bar{\sigma})\mu_{x-y}(\mbb{R}_+)\gamma_2(y)
+H\right]
\eeqlb
with $$
H=6k_0(l)+2\bar{\sigma}H_0\gamma_2(x_0).
$$
Now choose
\beqlb\label{conf of theta}
\theta=\max\left\{\frac{2(H+1)}
{\underline{\sigma}\kappa_0{\gamma_2(rx_0/2)}}+\frac{\bar{\sigma}}{\underline{\sigma}},
\frac{2(H+1)}{rx_0}\mathrm{e}^{c(l)l}\right\}.
\eeqlb
If $y\ge rx_0/2$, then $x-y\le (1-r/2)x_0$, and by \eqref{est of coup when jump (3)} and Condition \ref{cond noise} (B2-ii),
$$
\tilde{\mathcal{L}}_sF_{l}(x,y)
\le b_0(s)\Big(\frac{1}{2}\kappa_0(-\underline{\sigma}\theta+\bar{\sigma}) \gamma_2(rx_0/2)+H\Big).
$$
With the choice of $\theta$ in \eqref{conf of theta}, we have
$$
-\underline{\sigma}\theta+\bar{\sigma}\le-\frac{2(H+1)}{\kappa_0\gamma_2(rx_0/2)},
$$
which implies
$$
\tilde{\mathcal{L}}_sF_{l}(x,y) \le -b_0(s)
 \le-\frac{1}{2(\theta+1)}b_0(s)F_l(x,y).
$$
If $y< rx_0/2$, then $l>x-y> rx_0/2$. By \eqref{est of coup when jump (3)},
\beqlb\label{eq:4.33}
\tilde{\mathcal{L}}_sF_{l}(x,y)\le b_0(s)\left(-\phi(x)(x-y)\mathrm{e}^{-c(l)(x-y)}+H\right)
\le b_0(s)\left(-\frac{1}{2}\theta rx_0\mathrm{e}^{-c(l)l}+H\right).
\eeqlb
From \eqref{conf of theta},
$$
\theta\ge\frac{2(H+1)}{rx_0}\mathrm{e}^{c(l)l},
$$
so \eqref{eq:4.33} gives
$$
\tilde{\mathcal{L}}_sF_{l}(x,y)
\le-b_0(s)
\le-\frac{1}{2(\theta+1)}b_0(s)F_l(x,y).
$$

Combining the three cases completes the proof.\qed

Under Condition \ref{cond noise}, define
$$
G(x,y)=[V(x)+V(y)+\varepsilon F_{l}(x,y)]\mathbf{1}_{\{x\neq y\}}, \quad x,y\in\mbb{R}_+,
$$
where $l$ is sufficiently large such that $V(x)\ge(2\lambda_2+1)/\lambda_1$ for all $x\ge l$, with parameters $\lambda_1$ and $\lambda_2$ given in Condition \ref{cond noise} (C), $F_l$ is  given by \eqref{test function1} or \eqref{test function2}. The constant $\varepsilon$ equals $4\lambda_2/\lambda(l)$, where $\lambda(l)$ is the constant associated with $l$ from Proposition \ref{prop10} or Proposition \ref{prop11}, respectively.

Based on the preceding estimates, we now establish the following result.

\begin{proposition}\label{global est}
Suppose that Condition \ref{cond noise} holds. Then there exists $\bar{\lambda}>0$ such that for all $s\ge0$,
$$
\tilde{\mathcal{L}}_sG(x,y)\le-\bar{\lambda}b_0(s)G(x,y)\mathbf{1}_{\{x\neq y\}}.
$$
\end{proposition}
\proof Fix $s\ge0$. First consider the case of $x\neq y$ and $x \vee y \ge l$. By Propositions \ref{prop10} and \ref{prop11},
$\tilde{\mathcal{L}}_sF_l(x,y)\le0.$ Then, by Condition \ref{cond noise} (C),
\beqlb\label{final est1}
\tilde{\mathcal{L}}_sG(x,y)\ar=\ar
\tilde{\mathcal{L}}_sV(x)+\tilde{\mathcal{L}}_sV(y)
+\varepsilon\tilde{\mathcal{L}}_sF_l(x,y)\cr
\ar\le\ar
\left[-\lambda_1(V(x)+V(y))+2\lambda_2\right]b_0(s)\cr
\ar\le\ar
\left[-\frac{1}{2}\lambda_1(V(x)+V(y))-1\right]b_0(s)\cr
\ar\le\ar
\left[-\frac{1}{2}\lambda_1(V(x)+V(y))-(2(\theta+1)\varepsilon)^{-1}\varepsilon F_l(x,y)\right]b_0(s)\cr
\ar\le\ar
-\frac{1}{2}(\lambda_1\wedge(\varepsilon(\theta+1))^{-1})b_0(s)G(x,y).
\eeqlb
Next, consider $x \vee y \le l$ and $x\neq y$. Then $0 < |x-y|\le l$. From Propositions \ref{prop10} and \ref{prop11},
$$
\tilde{\mathcal{L}}_sF_l(x,y)\le-\lambda(l)b_0(s)F_l(x,y).
$$
Hence, by Condition \ref{cond noise} (C),
\beqnn
\tilde{\mathcal{L}}_sG(x,y)\ar=\ar\tilde{\mathcal{L}}_sV(x)+\tilde{\mathcal{L}}_sV(y)+\varepsilon\tilde{\mathcal{L}}_sF_l(x,y)\cr
\ar\le\ar
\left[-\lambda_1(V(x)+V(y))+2\lambda_2-\varepsilon\lambda(l) F_l(x,y)\right]b_0(s).
\eeqnn
Note that $F_l>1$ (since $\phi > 1$ and $\psi_l \ge1$), so
$$
2\lambda_2-\varepsilon\lambda(l)F_l(x,y)/2
\le2\lambda_2-\varepsilon\lambda(l)/2=0.
$$
This implies
\beqlb\label{final est2}
\tilde{\mathcal{L}}_sG(x,y)\ar\le\ar\left[-\lambda_1(V(x)+V(y))
-\frac{1}{2}\varepsilon\lambda(l)F_l(x,y)
-2^{-1}\varepsilon\lambda(l)F_l(x,y)+2\lambda_2\right]b_0(s)\cr
\ar\le\ar
\left[-\lambda_1(V(x)+V(y))-\frac{1}{2}\varepsilon\lambda(l)F_l(x,y)\right]b_0(s)\cr
\ar\le\ar-\frac{1}{2}[\lambda_1\wedge(\varepsilon\lambda(l))]G(x,y)b_0(s).
\eeqlb

From \eqref{final est1} and \eqref{final est2}, we conclude that
$$
\tilde{\mathcal{L}}_sG(x,y)\le-\bar{\lambda}b_0(s)G(x,y)\mathbf{1}_{\{x\neq y\}}
$$
with
$$
\bar{\lambda}=\frac{1}{2}\min \left\{\lambda_1, \varepsilon\lambda(l), \frac{1}{\varepsilon(\theta+1)}\right\}.
$$
The proof is complete.\qed

We now prove Theorem \ref{exponential ergod}.

\noindent{\it  Proof of Theorem \ref{exponential ergod}}. Fix $s\ge0$ and  let $\mathbf{E}^{(X_s,Y_s)}(\cdot):=\mathbf{E}(\cdot|\sigma(X_s,Y_s))$. For $n\ge1$, define
$$
T_{n,s}=\inf\left\{t\ge s: |X_t-Y_t|\notin[1/n,n]\right\}.
$$
Since the  coupled process $(X,Y)$ is non-explosive, we have $T_{n,s}\uparrow T_s$ almost surely as $n\rightarrow\infty$, where
$$
T_s=\inf\left\{t\ge s: X_t=Y_t\right\}
$$
is the coupling time starting from time $s$. For $t\ge s$, integration by parts yields
\beqnn
\ar\ar\mathbf{E}^{(X_s,Y_s)}\left[\mathrm{e}^{\bar{\lambda}\int_{s}^{t\wedge T_{n,s}}b_0(r)\,\dd r}G(X_{t\wedge T_{n,s}},Y_{t\wedge T_{n,s}})\right]\cr
\ar=\ar G(X_s,Y_s)+\mathbf{E}^{(X_s,Y_s)}\left(\int_s^{t\wedge T_{n,s}}\mathrm{e}^{\bar{\lambda}\int_{s}^{r\wedge T_{n,s}}b_0(r)\,\dd r}\left[\bar{\lambda}b_0(r)G(X_r,Y_r)+\tilde{\mathcal{L}}_rG(X_r,Y_r)\right]\,\dd r\right)\cr
\ar\le\ar G(X_s,Y_s),
\eeqnn
where $\bar{\lambda}$ is given in Proposition \ref{global est}, and the inequality follows from Proposition \ref{global est}. (Similar arguments appear in \cite[Theorem 3.1]{LW19}.)  By Fatou's lemma, letting $n\rightarrow\infty$ gives
$$
\mathbf{E}^{(X_s,Y_s)}\left[\mathrm{e}^{\bar{\lambda}\int_{s}^{t\wedge T_{s}}b_0(r)\,\dd r}G(X_{t\wedge T_{s}},Y_{t\wedge T_{s}})\right]\le G(X_s,Y_s).
$$
By our convention that $X_t=Y_t$ for $t\ge T_{s}$, we have $G(X_t, Y_t)=0$ for all $t\ge T_s$, which implies
\beqnn
\mathbf{E}^{(X_s,Y_s)}\left[\mathrm{e}^{\bar{\lambda}\int_{s}^{t\wedge T_{s}}b_0(r)\,\dd r}G(X_{t\wedge T_{s}},Y_{t\wedge T_{s}})\right]
\ar=\ar \mathrm{e}^{\bar{\lambda}\int_{s}^{t}b_0(r)\,\dd r}\mathbf{E}^{(X_s,Y_s)}\left[G(X_{t},Y_{t})\mathbf{1}_{\{T_s>t\}}\right]\cr
\ar=\ar
\mathrm{e}^{\bar{\lambda}\int_{s}^{t}b_0(r)\,\dd r}\mathbf{E}^{(X_s,Y_s)}\left[G(X_{t},Y_{t})\right].
\eeqnn
Therefore,
\beqlb\label{a coupling est of XY}
\mathbf{E}^{(X_s,Y_s)}\left[G(X_{t},Y_{t})\right]\le\mathrm{e}^{-\bar{\lambda}\int_{s}^{t}b_0(r)\,\dd r}G(X_s,Y_s).
\eeqlb
Note that $0<\inf_{x\neq y}F_l(x,y)\le \sup_{x\neq y}F_l(x,y)<\infty$. From the definition of $d_V$ in \eqref{metric of dv}, there exist constants $c_1>c_2>0$ such that
$$
c_1d_V(x,y)\le G(x,y)\le c_2d_V(x,y).
$$
Hence, by \eqref{wasserstein dist} and \eqref{a coupling est of XY},
\beqnn
W_V(P_{s,t}(x,\cdot),P_{s,t}(y,\cdot))\ar\le\ar c_1^{-1}\inf_{\pi\in{\mathcal{C}}(P_{s,t}(x,\cdot),P_{s,t}(y,\cdot))}
\int_{\mathbb{R}_+^2}G(z,\tilde{z})\pi(\dd z,\dd \tilde{z})\cr
\ar\le\ar c_1^{-1}\mathrm{exp}\left(-\bar{\lambda}\int_{s}^{t}b_0(r)\,\dd r\right)G(x,y)\cr
\ar\le\ar c_2 c_1^{-1} \mathrm{exp}\left(-\bar{\lambda}\int_{s}^{t}b_0(r)\,\dd r\right)d_V(x,y).
\eeqnn
The proof is complete.\qed

\section{Appendix}

In this section, we proof Proposition \ref{ii1.9}, which adapts the Yamada-Watanabe argument to our setting.

\noindent{\it Proof of Proposition \ref{ii1.9}}. Consider a sequence of functions $\{\phi_k\}$ on $\mathbb{R}$ as follows. For each $k\ge0$, define $a_k=\exp\{-k(k+1)/2\}$. Then $a_k\rightarrow0$ decreasingly as $k\rightarrow\infty$ and $\int_{a_k}^{a_{k-1}}z^{-1}\,\dd z=k$ for $k\ge1$. Let $x\rightarrow\psi_k(x)$ be a nonnegative continuous function on $\mathbb{R}$ with support in $(a_k,a_{k-1})$ and satisfying  $\int_{a_k}^{a_{k-1}}\psi_k(x)\,\dd x=1$ and $0\le\psi_k(x)\le2(kx)^{-1}$ for $a_k<x<a_{k-1}$. For each $k\ge1$, define
$$
\phi_k(z)=\int_0^{|z|}\dd y\int_0^y\psi_k(x)\dd x,\quad z\in\mathbb{R}.
$$
The sequence $\{\phi_k\}$ satisfies the following properties:

(a) $\phi_k(z)\rightarrow|z|$ non-decreasingly as $k\rightarrow\infty$;

(b) $0\le \phi'_k(z)\le1$ for $z\ge0$ and $-1\le\phi'_k(z)\le0$ for $z\le0$;

(c) $0\le |z|\phi''_k(z)=|z|\psi_k(|z|)\le2k^{-1}$ for $z\in\mathbb{R}$.

For $xz\ge0$, we have
\beqlb\label{est of phi_k2}
D_z\phi_k(x)\le z,
\eeqlb
and by Taylor expansion,
\beqlb\label{est of phi_k3}
\begin{array}{ll}
\displaystyle
D_z\phi_k(x)
=z^2\int_0^1\psi_k(x+tz)(1-t)\dd t
\le z^2\int_0^1\frac{2}{k(x+tz)}(1-t)\dd t\le\frac{z^2}{kx}.
\end{array}
\eeqlb
Suppose $X$ and $\tilde{X}$ are two solutions of \eqref{ii1.5} with $X_0=\tilde{X}_0$. Define $\Upsilon_t=X_t-\tilde{X}_t$ for $t\ge0$. From \eqref{ii1.5}, we obtain
\beqlb\label{variate of X}
\Upsilon_t\ar=\ar\int_0^t\left[\gamma_0(s-,X_{s-})-\gamma_0(s-,\tilde{X}_{s-})\right]\dd s
+\int_0^t\int_{\gamma_1(\tilde{X}_{s-})}^{\gamma_1(X_{s-})}\,W(\dd s,\dd u)\cr
\ar\ar\,\,
+\int_0^t\int_0^\infty\int_{\gamma_2(\tilde{X}_{s-})}^{\gamma_2(X_{s-})}z\mathbf{1}_{\{\Upsilon_{s-}>0\}}\,\tilde{N}(\dd s,\dd z,\dd u)\cr
\ar\ar\,\,
-\int_0^t\int_0^\infty\int_{\gamma_2(X_{s-})}^{\gamma_2(\tilde{X}_{s-})}z\mathbf{1}_{\{\Upsilon_{s-}\le0\}}\,\tilde{N}(\dd s,\dd z,\dd u).
\eeqlb
 Fix $m\ge1$. Recall the stopping times
\beqnn
\tau^-_{1/m}=\inf\{t\ge0:X_t\ge 1/m\},\quad
\tau^+_{m}=\inf\{t\ge0:X_t\le m\},\quad \tau_m=\tau^-_{1/m}\wedge\tau^+_m.
\eeqnn
Define $\tilde{\tau}_m$ analogously for $\tilde{X}$ and set $\sigma_m=\tau_m\wedge\tilde{\tau}_m$. Applying It\^{o}'s formula to \eqref{variate of X} yields
\beqnn
\phi_k(\Upsilon_{t\wedge\sigma_{m}})
\ar=\ar \int_0^{t\wedge\sigma_{m}}b_1(s)\phi''_k(\Upsilon_{s-})
\left|\gamma_1(X_{s-})-\gamma_1(\tilde{X}_{s-})\right| \dd s\cr
\ar\ar\,
+\int_0^{t\wedge\sigma_{m}}\phi'_k(\Upsilon_{s-})\left[\gamma_0(s-,X_{s-})-\gamma_0(s-,\tilde{X}_{s-})\right]\dd s\cr
\ar\ar\,
+\int_0^{t\wedge\sigma_{m}}\int_0^\infty\int_{\gamma_2(\tilde{X}_{s-})}^{\gamma_2(X_{s-})}
D_z\phi_k(\Upsilon_{s-})\mathbf{1}_{\{\Upsilon_{s-}>0\}}b_2(s)\,\dd s\mu(\dd z)\dd u\cr
\ar\ar\,
+\int_0^{t\wedge\sigma_{m}}\int_0^\infty\int_{\gamma_2(X_{s-})}^{\gamma_2(\tilde{X}_{s-})}
D_{-z}\phi_k(\Upsilon_{s-})\mathbf{1}_{\{\Upsilon_{s-}\le0\}}b_2(s)\,\dd s\mu(\dd z)\dd u+mart..
\eeqnn
Using the properties of $\phi_k$, \eqref{est of phi_k3} and assumptions on the functions $\gamma_0,\gamma_1$ and $\gamma_2$, we have
\beqnn
\phi_k(\Upsilon_{t\wedge\sigma_{m}})
\ar\le\ar K\left\{\frac{1}{k}\int_0^{t\wedge\sigma_{m}}
\left(2b_1(s)+\left(\int_0^1z^2\,\mu(\dd z)\right)b_2(s)\right)\dd s+\int_0^{t\wedge\sigma_{m}}b_0(s-)|\Upsilon_{s-}|\dd s\right.\cr
\ar\ar\,\,
\left.+\left(\int_1^\infty z\,\mu(\dd z)\right)\int_0^{t\wedge\sigma_{m}}b_2(s)|\Upsilon_{s-}|\,\dd s\right\}+mart.
\eeqnn
for some $K>0$. Taking expectations and letting $k\rightarrow\infty$ gives
$$
\mathbf{E}\left[\left|\Upsilon_{t\wedge\sigma_{m}}\right|\right]\le K\int_0^t(b_0(s)+b_2(s))\mathbf{E}[|\Upsilon_{s-\wedge\sigma_{m}}|]\,\dd s.
$$
By Gronwall's inequality,
$$
\mathbf{E}[|\Upsilon_{t\wedge\sigma_{m}}|]=0
$$
holds for all $t\ge0$. Since $\sigma_{m}\rightarrow\tau\wedge\tilde{\tau}$ as $m\rightarrow\infty$, the right continuity of $t\mapsto\Upsilon_t$ implies the desired result.
\qed

%%%%%%%%%%%%%%%%%%%%%%%%%%%%%%%%%%%%%%%%%%%%%%%%%%%%%%%%%%%%%%%%%%%
%%                                                               %%
%% Use the two commands below for producing your bibliography    %%
%% with bibtex, then comment again the commands and include the  %%
%% content of the .bbl file in this file below the commands.     %%
%%                                                               %%
%%%%%%%%%%%%%%%%%%%%%%%%%%%%%%%%%%%%%%%%%%%%%%%%%%%%%%%%%%%%%%%%%%%

%\bibliographystyle{amsplain}
%\bibliography{yourbibfilename}

% add below the content of your .bbl file produced by bibtex.

{\bf   {Acknowledgements.}}
The research of Shukai Chen is supported
    by the National Key R\&D Program of China
(No.2022YFA1006003), NSFC grant of China 
(No.12401167), 
Fujian Provincial Natural Science Foundation of China (No.2024J08050) and the Education and Scientific Research Project for Young and Middle-aged Teachers in Fujian Province of China 
(No.JAT231015). The research of Xu Yang is supported by NSFC grant of China 
(No.12471135). The research of Xiaowen Zhou is supported by NSERC (RGPIN-2021-04100).

\end{document}